\documentclass[12pt]{amsart}
\usepackage{amsmath,amssymb,amsthm,
color, euscript, enumerate}
\usepackage{dsfont}
\usepackage{longtable}
\newcommand{\RR}{{\mathbb R}}

\newcommand{\be}{\begin{equation}}
\newcommand{\ee}{\end{equation}}
\newcommand{\ba}{\begin{array}}
\newcommand{\ea}{\end{array}}
\newcommand{\bea}{\begin{eqnarray}}
\newcommand{\eea}{\end{eqnarray}}

\usepackage{graphicx}
\usepackage{color}
\newcommand{\corr}[1]{{\color{black}{#1}}}
\usepackage{transparent} 

\newtheorem{theorem}{Theorem}  

\newtheorem{lemma}{Lemma}[section]

\DeclareMathOperator{\res}{res}

\numberwithin{equation}{section}
\allowdisplaybreaks

\begin{document}

\begin{center}
{\large   \bf Stability of KdV solitons  } 
 
\vskip 15pt

{ {\large Derchyi Wu}\\
Institute of Mathematics, Academia Sinica, 
Taipei, Taiwan}

\vskip 10pt

{\today}

\vskip 10pt
{\bf Abstract}
\end{center}
\begin{enumerate}
\item[]{\small We prove an orbital stability theorem of KdV $n$-solitons with explicit phase shifts in the soliton region with cones around the $x$-axis and lines determined by bound states of the KdV $n$-solitons removed. }
\end{enumerate} 
%%%%%%%%%%%%%%%%%%%%%%%%%%%%%%%%%%%%%%%%%%%%%%%%%%%%%%%%%%%%%%%%%%%%%%%%
\section{Introduction}\label{S:introduction}
%%%%%%%%%%%%%%%%%%%%%%%%%%%%%%%%%%%%%%%%%%%%%%%%%%%%%%%%%%%%%%%%%%%%%%%%

%The  Korteweg-de Vries ({\bf KdV}) equation derived by Korteweg and de Vries \cite{KdV95} in 1895, is  an asymptotic model for the propagation of one-dimensional small amplitude  long  surface waves. Special solitary solutions which represent  a balance of dispersive and non-linear effects of water waves  were first observed by Russel in 1834 \cite{Mi76} and computed by Boussinesq \cite{Bou71,Bou72} and Rayleigh independently \cite{Ray76}.   Modern KdV theory as well as the integrable system theory were unleashed by the discovery that the KdV equation can be linearised by the inverse scattering theory ({\bf IST}) of the Schr$\ddot{\mbox{o}}$dinger operator \cite{ZK65, Zab67, Lax68}
%\[-f'' +u f   =k^2f  \]in 1960's. Here $'$, $''$ denote $\partial/\partial x$, $\partial^2/\partial x^2$, and $k\in\CC$. Famous rigorous (and beautiful!) IST theories of the KdV equation are proved by Faddeev \cite{F59, F64, F76} using the Gelfand-Levitan-Marchenko equation and by Deift and Trubowitz via Backlund transformation and a nonlinear first order system in terms of the trace formula \cite{DT79}.

The solitary solutions 
\be\label{E:line-tau-oblique}
 u (x,t) 
= -\frac {\beta ^2 }2\textrm{sech}^2\frac 12(\beta  x-4\beta ^3 t +\ln \alpha )
\ee  are the simplest examples in  an important class of  solutions of the Korteweg-de Vries ({\bf KdV}) equation 
\be\label{E:KdV}
 u_t-6uu_{x}+u_{xxx}=0,  
\ee called the KdV $n$-solitons. They  correspond to reflectionless KdV solutions in the inverse scattering theory ({\bf IST}) \cite{F59, F64, F76, DT79} and possess the following explicit form 
  \be\label{E:n-soliton-tau}
u_{0}(x,t)=-2\frac{d^2}{dx^2}\log \textrm{Wr}(a_1(x,t),\cdots,a_n(x,t))
\ee with               
\begin{gather}
\textrm{Wr} (a_1(x,t),\cdots,a_n(x,t))=\left|
\begin{array}{llll}
a_{1} &a_{2} & \cdots & a_{n}\\
a'_{1} &a'_{2} & \cdots & a'_{n}\\
\vdots & \vdots &\ddots &\vdots\\
a^{(n-1)}_{1} &a^{(n-1)}_{2} & \cdots & a^{(n-1)}_{n}
\end{array}
 \right|,\label{E:n-soliton}\\
 a_j(x,t)=(-1)^{j+1}e^{-\beta_j (x-4\beta_j^2 t)}+\alpha_je^{+\beta_j (x-4\beta_j^2 t)},\label{E:n-soliton-a_j} \\
\alpha_j>0, \ 0<\beta_1<\cdots<\beta_n,\label{E:n-soliton-alpha}
\end{gather}and $a'_j=\frac{\partial }{\partial x }a_j,\,a''_j=\frac{\partial^2 }{\partial x^2 }a_j,\,a^{(h)}_j=\frac{\partial^h}{\partial x^h}a_j,\,  1\le j\le n$.  They represent non linear superpositions of $n$ travelling waves with different speeds $\beta_j$ which interact elastically   and behave asymptotically as the sum of $n$ travelling waves as $t\to +\infty$ \cite{KM56, GGKM67, H71,Z71,T72,WT72,GGKM74,DT79}.

Via variational methods of Hamiltonian systems,  asymptotic characterization, concentration compactness techniques, and nonlinear methods, it is shown  orbital  stability in various Sobolev spaces is vaild \cite{Ben72,Bo75,Wein86,MS93,MMT02,MV03,MM05,AMV13,KV22}. Namely, they proved if
\[
\Vert u_v(x,0)-u_0(x,0)\Vert <\epsilon \]where $\Vert\cdot\Vert$ represents  various Sobolev types norms ($H^1(\RR)$ in \cite{Ben72} and  $H^{-1}(\RR)$ in \cite{KV22} for instance), then for some constant $C$,
 \be\label{E:KdV-orbital}
\begin{gathered}
  \sup_{t> 0,x\in\RR}\ \inf_{x_j(t)\in\RR}\Vert u_v(x,t)+2\frac{d^2}{dx^2}\log \textrm{Wr}\left(a_1(x-x_1(t),t),\cdots,a_n(x-x_n(t),t)\right)\Vert<C\epsilon.
\end{gathered}
\ee  The necessity of phase shifts $x_j(t)$   in \eqref{E:KdV-orbital} can be seen from the famous counterexample   
     \[ 
     \begin{gathered}
         u_j (x,t) 
=  -2\beta_j ^2 \textrm{sech}^2(\beta_j (x-4\beta_j ^2t)),\ \ \beta_1\ne\beta_2,\ \ 
|\beta_1-\beta_2|\ll 1.
     \end{gathered}
\] Hence $ | u_1(x,0)-u_2(x,0)|_{L^\infty}\ll 1$ but $
     |u_1 (x,t) -u_2 (x,t)|_{L^\infty}\sim\mathcal O(1)$  as $t=\frac 1{\beta_2-\beta_1}$.   {\color{black}Detailed description of the phase shifts 
      is not available   except taking $t\to\infty$   \cite{M05, MMT02, KP05, MT14}. %It is also important to compare with the asymptotic theorem by taking $t\to\infty$ along lines of fixed slopes  for any fixed KdV solution via IST approach \cite{ES83}.} 

 The goal in this paper is to provide an orbital stability with explicit phase shifts of the KdV $n$-solitons in a regime where $t$ is not necessarily large. The inverse scattering theory approach adopted in this paper also exhibits the inner working in a different and clearer light.

\begin{theorem}\label{T:stability} Given a KdV $n$-soliton $u_{0}(x,t)$, positive constants   $\mathfrak a$, $ \sigma>2$,   $  \color{black}\tau $,  
 there exists  $C\gg 1$   such that if   the   KdV solution $u_v(x,t)  $ satisfies 
 \[\begin{array}{c}
u_v(x,0)=u_0(x,0 )+v(x), \quad 
\sum_{h=0}^2|v^{(h)} e^{\mathfrak a|x|}|_{ L^\infty }\le     \frac{\epsilon^\sigma}C     
\end{array}\]for $\color{black}\epsilon < \min_{1\le j\ne j'\le n}\{\mathfrak a, \frac {\beta_{0,1}}2,\frac {|\beta_{0,j}-\beta_{0,j'}|}2 \} $, then    
\begin{multline*}
   |u_v(x,t) 
   +2\frac{d^2}{dx^2}\log \textrm{Wr} (a_1(x-x_1(t),t),\cdots\\\hskip1in\cdots,a_n(x-x_n(t),t) )|_{  L^\infty(\corr{ \frac x{4t}\in \cap_{j=1}^n \{c^2:|c-\beta_{0,j}|\ge \tau,  \frac 1\tau\ge c\ge\epsilon\},t>0 })}\\ < C\corr{\epsilon^{\min\{1,\sigma-2\}}}  
\end{multline*}
with
\[
\begin{split}
x_{j}(t)=&4(\beta_{v,\mu(v)-n+j}^2-\beta_{0,j}^2)t,\quad j=1,\cdots,n,\quad n\le \mu(v),\\
a_j(x,t)=&(-1)^{j+1}e^{-\beta_{0,j}(x-4\beta_{0,j}^2t)}+\alpha_{0,j} e^{+\beta_{0,j}(x-4\beta_{0,j}^2t)},\\
\alpha_{0,j}=&\frac{2(-1)^{j-1}\beta_{0,j}}{\gamma_{0,j}}\Pi_{l\ne j}\frac{\beta_{0,l}+\beta_{0,j}}{\beta_{0,l}-\beta_{0,j}}
\end{split}
\] and
       $\beta_{v,j}$, $\gamma_{v,j}$    are bound states, norming constants  of $u_v(x,0)$.%, and $\beta_1,\cdots,\beta_{v,\mu(v)}$   are bound states   of $u_v(x,0)$. 
\end{theorem}

For simplicity, throughout the paper, $C$ denotes a uniform constant  which is independent of $x$, $t$,  and the spectral variable $k$, and  $C_{x_+}=C_{x_+}(x) $ ($C_{x_-}=C_{x_-}(x) $) is a non increasing (decreasing) function.  

 We remark that IST approach has been applied to  justify a global $L^\infty$-stability  with no phase shift of KP multi-line solitons recently \cite{Wu22}, and an $L^\infty$-orbital stability with explicit phase shifts,  in a regime uniform with respect to $t$, of nonlinear Schr$\ddot{\mbox o}$dinger (NLS) $1$-solitons \cite{CP14}. %We note that $L^\infty$ For the NLS, \cite{CP14} adopted the deepest descent method and Backlund transformation  theory.

Besides, our result differs from the immense and descent literatures about KdV asymptotics where $t\to\infty$ is demanded (see \cite{ES83, DVZ94, GT09, AMSM23} for instance). Major distinction and technical difficulty are due to,  under arbitrarily small perturbation, the number of bound states could increase   \corr{and there could be no uniform positive lower bound for the bound states $\beta_{v,j}$} \cite{D83}. 
To remove these  difficulties, we shall adopt a deformed {\bf Gelfand-Levitan-Marchenko (GLM)} equation method, enabled by the exponentially decaying property of the perturbation and   introduced by Eckhaus-Schuur \cite{ES83}, \corr{and only consider stability in the soliton region with cones around the $x$-axis and lines determined by bound states of $u_0(x,0)$ removed \cite{DVZ94}}. 

%Unlike the NLS $1$-soliton, and also the immense and descent literatures about KdV asymptotics \cite{???}, major obstructions to deriving the KdV $n$-soliton stability via IST include lack of   global (in $t\ge 0$) $L^1_2(\RR^\pm)$-estimates of Fourier transform   of the reflection coefficients $R_\pm$ (see \eqref{E:cohn} or \cite{C79,C82})   even for a single fixed initial data;    

The paper is organized as follows:  we   give a brief review of the IST of the KdV equation for convenience in Section \ref{S:KdV-IST}.   
In Section \ref{S:perturbation}, for the direct scattering problem, we shall characterize properties and derive estimates of the scattering data under the perturbation. For the inverse scattering problem,  we shall formulate a deformed GLM equation with a kernel which removes the the smaller bound states {\color{black}in the soliton region $x\ge 4\epsilon^2 t$.

In Section \ref{S:stability}, we prove Theorem \ref{T:stability} by analyzing more about the deformed GLM equation. More precisely, we  decompose the solution of the GLM equation into the continuous part $B_{v,c}$ and the discrete part $B_{v,d}$. A refinement of the GLM equation of KdV $n$-solitons \cite[Lemma 6.1]{ES83} shows that uniform estimates about $B_{v,c}$ hold only in the soliton region outside cones around the $x$-axis and lines determined by bound states of  $u_0(x,0)$. Estimates about  $B_{v,d}$ are done by using the Wronskian represenation formula  \eqref{E:n-soliton} and an application of the asymptotic theory of KdV $n$-solitons \cite{K04, K18}.  }

{\bf Acknowledgments}.  We are indebted to inputs from J. H. Chang,  K. \corr{Gkogkou}, Y. Kodama, and P. Miller.  We would like to give special thanks to A. Calini  and the hospitality of Loughborough University and Isaac Newton Institute, UK in September 2022.  
This research   was   supported by NSC 111-2115-M-001 -003 -. 

%%%%%%%%%%%%%%%%%%%%%%%%%%%%%%%%%%%%%%%%%%%%%%%%%%%%%%%%%%%%%%%%%%%%%%%%%
\section{The IST for the KdV equation}\label{S:KdV-IST}
%%%%%%%%%%%%%%%%%%%%%%%%%%%%%%%%%%%%%%%%%%%%%%%%%%%%%%%%%%%%%%%%%%%%%%%%%

We give a brief review of the IST of KdV equation \cite{F59, F64, DT79} in this section. 

For the direct problem, suppose $u\in L^1_2=\{\varphi: |(1+|x|^2)\varphi |_{L^1}<\infty\}$, for $k\in \mathbf R\backslash\{0\}$, there exist  eigenfunctions of  the Schr$\ddot{\mbox{o}}$dinger equation
\begin{align}
&-f''_\pm +u (x) f_\pm =k^2f_\pm,\label{E:eigenfunction} \\
&f_{\pm}(x,k)=m_{\pm}(x,k)e^{\pm ikx} \quad\to \quad e^{\pm ikx}      \textit{ as $x\to \pm\infty$}, \nonumber 
\end{align} 
such that   $m_\pm(x,\cdot)-1 =  \pm\int_0^{\pm\infty}B_\pm(x,y)e^{\pm2iky}dy $  are in the Hardy space $H^{2+} $,    
\begin{align}
|  m_{\pm}(x,k)  -1 | 
\le  &    \frac{C(1+\max(\mp x,0))}{1+|k| }| u|_{L^1_2},\label{E:m-infty}\\
 |B_\pm(x,\cdot) |_{L^\infty(\mathbf R^\pm)\cap L^1(\mathbf R^\pm)}\le & C_{x_\pm}| u|_{L^1_2}. \label{E:fourier-m}
\end{align}

Hence we obtain  $\textbf {transmission}$ and $\textbf {reflection}$ coefficients $T_+$, $R_+$
\be\label{E:RT-coefficient} 
f_{\mp}(x,k)= \frac{R_\pm (k)}{T_\pm (k)}f_{\pm}(x,k)+\frac{1}{T_\pm (k)}f_{ \pm}(x,-k) ,\ k\in\mathbf R\backslash \{0\},
\ee and define  the {\bf forward scattering transform} by 
\be\label{E:forward}  
\mathcal S(u(x,0))\equiv \{R_+(k), i\beta_{ 1},\cdots,i\beta_{n},\gamma_{+, 1},\cdots,\gamma_{+,n} \}, 
\ee where    
   $
i\beta_j$  are simple poles of $T_+(k)$, $\beta_j>0 $, and $
\gamma_{\pm,j}= \left(\int_{-\infty}^\infty f_\pm ^2(x,i\beta_j)dx \right)^{-1}$. In particular, for the $n$-soliton $u _0(x,0)$ defined by \eqref{E:n-soliton-tau}-\eqref{E:n-soliton-alpha}, one has \cite[Theorem 3.6]{DT79}
   \be \label{E:pole-soliton}
  \begin{split}
  &R_+\equiv 0, \quad
   \gamma_{ j}= \frac{2(-1)^{j-1}\beta_j}{\alpha_j}\left(\Pi_{1\le l\le   n,l\ne j}\frac{\beta_l+\beta_j}{\beta_l-\beta_j}\right),\quad 1\le j\le   n. 
  \end{split}
  \ee  Moreover, 
the algebraic and analytic constraints hold,
\begin{equation} \label{E:alg-con}
\begin{split}
R_\pm(k)=\overline{R_\pm(-k)} ,&\quad T_\pm(k)=\overline{T_\pm(-k)}=T_\mp(k)\equiv T(k) ,\\ 
|T(k)|^2+|R_\pm(k)|^2=1,&\quad R_+(k)T(-k)+R_-(-k)T(k)=0,\\
&\quad\gamma_{+,j}\gamma_{-,j}=-\left[ (1/T  )'\right]^{-2}(i\beta_j),
\end{split}
\end{equation}and  
\be\label{E:analytic-con}\begin{gathered}
 (T-1)^{-1}\in  H^{2+} ,\       |k ^{N +1}R_\pm(k) |_{L^\infty} \le      C\sum_{h=0}^N|u^{(h)}|_{L^1_2} ,   \\
        \pm  \int_x^{\pm\infty}(1+s^2)|\frac{d}{ds}F_\pm | ds,\ \pm \int_x^{\pm\infty} |F_\pm |ds       \le    C_{x_\pm}|u|_{L^1_2},    
\end{gathered} \ee
with 
\be\label{E:GLM-kernel-0}
 F_\pm(x)= \frac 1\pi\int_{-\infty}^\infty R_\pm(k) e^{\pm 2ikx}dk.   
 \ee 
%The forward scattering transform \eqref{E:forward} is not continuous in general. Since the number of the bound states could increase under a perturbation \cite{D83}. We shall discuss perturbation properties of $\mathcal S$ in Section \ref{S:perturbation}.

We turn to  the inverse problem. Given a scattering data
$
\mathcal S = \{R_+(k),i\beta_{ 1},\cdots,i\beta_{n}$, $\gamma_{+, 1},\cdots,\gamma_{+,n} \}$ satisfying % and the algebraic and analytical constraints
\begin{equation} \label{E:inv-sd}
\begin{gathered}
0<\beta_1<\cdots<\beta_n,\quad 0<\gamma_{+, 1},\cdots,\gamma_{+,n} ,\\
R_+(k)=\overline{R_+(-k)} ,\quad 0<|R_+|<1,\quad  |R_+|\le \frac{C}{1+|k|^{ 2}} , \\
  \int_x^{+\infty}(1+s^2)|\frac{d}{ds}F_+ | ds,\quad  \int_x^{+\infty} |F_+ |ds   \le    C_{x_+} ,
\end{gathered}
\end{equation} Define  $R_-$, $\gamma _{-,j}$, and $F_-$ by \eqref{E:alg-con}, \eqref{E:GLM-kernel-0}, and 
   \begin{equation}   \label{E:GLM-kernel}
\Omega_\pm(x)=  F_\pm(x)+2\sum_{j=1}^n\gamma_{\pm,j}\exp(\mp 2\beta_{j}x).  
\end{equation} Then the {\bf  Gelfand-Levitan-Marchenko} equations
  \begin{equation}  \label{E:GLM}
\Omega_\pm(x+y)+B_\pm(x,y)\pm\int_0^{\pm\infty} \Omega_\pm (x+y+u)B_\pm(x,u)du=0,\  y\in \mathbf R^\pm 
\end{equation}admit solutions $B_\pm(x,y) $ satisfying
\[ 
\begin{split}
    |B_\pm(x,\cdot)|_{L^\infty(\mathbf R^\pm)\cap L^1(\mathbf R^\pm)}\le & C_{x_\pm} C_{\gamma_j,\beta_j} 
\end{split}
\] where $C_{\gamma_j,\beta_j} $ a constant determined by $\gamma_j,\beta_j$.
Define 
 $ 
 m_\pm=1  \pm\int_0^{\pm\infty}B_\pm(x,y)e^{\pm2iky}dy $ and $f_\pm(x,k)= m_{ \pm }(x,k)e^{\pm ikx}$. Then  
\be\label{E:inverse} 
\begin{gathered}
      - f'' _\pm   +u(x) f_\pm\ =   \ k^2f_\pm  ,\\
      u(x)=\mp\partial B_\pm(x,0^\pm)/\partial x,
\end{gathered}\ee with\be\label{E:inverse-u}
 |u|_{L^1_2}\le  C_{\gamma_j,\beta_j} .\ee

  The {\bf inverse scattering transform} is defined by 
  \be\label{E:inverse-scattering-transform} u(x)=\mathcal S^{-1}( \{R_+(k),i\beta_{ 1},\cdots,i\beta_{n},\gamma_{+, 1},\cdots,\gamma_{+,n} \})\equiv \mp\partial B_\pm(x,0^\pm)/\partial x.\ee
%Moreover, equipped with proper topology, the inverse scattering transform $\mathcal S^{-1}$  is continuous.
 
The linearization theorem of the KdV equation is, suppose  
 \begin{equation}\label{E:evolution-kdv}
 \begin{gathered}
  %- f'' _\pm   +u(x,t) f_\pm  =  \ k^2f_\pm,\\
  u_{t}-6uu_{x }+u_{xxx} =0,\quad u(x,0)\in L^1_2,\\
  \mathcal S(u(x,t))\equiv \{R_+(k,t),i\beta_{ 1}(t),\cdots,i\beta_{n}(t),\gamma_{+, 1}(t),\cdots,\gamma_{+,n} (t)\},
  \end{gathered}\ee
    one has   
\begin{equation}\label{E:evolution}
\begin{array}{rcl}
{  R_{\pm}(k,t)}= & {R_\pm (k,0)}\exp(\pm 8ik^3t)&{  \equiv {R_\pm (k)}\exp(\pm 8ik^3t) ,}\\
{  \beta_j (t)}= & \beta_j(0)&{  \equiv \beta_j,} \\
{  \gamma_{\pm,j}(t)}= & \gamma_{\pm,j}(0) \exp(\pm 8\beta_j ^3t)&{  \equiv \gamma_{\pm,j}  \exp( { \pm 8\beta_j ^3t})
} 
\end{array}
\end{equation} \cite{GGKM67,T74}. 
As a result, let $F_\pm(x,t)=\widehat{R_{\pm}(\cdot,t)}$ and 
             \begin{equation} \label{E:GLM-kernel-t}
\Omega_\pm(x,t)=  F_\pm(x,t)+2\sum_{j=1}^n\gamma_{\pm,j}\exp(\mp 2\beta_{j}x{ \pm 8\beta_j^3 t})   . 
\end{equation} 
  Since  
             \[
                \int_x^{\pm\infty}(1+|s|^2)\left|\frac{\partial F_\pm}{\partial s}(s,t)\right|ds\le c_{\pm}(t)C_{x_\pm} .  
              \]  Applying the above IST, one can solve the GLM equation 
               \begin{equation}  \label{E:GLM-t}
\Omega_\pm(x+y,t)+B_\pm(x,y,t)\pm\int_0^{\pm\infty} \Omega_\pm (x+y+u,t)B_\pm(x,u,t)du=0.\  y\in \mathbf R^\pm 
\end{equation}and Cauchy problem of the KdV equation with $L^1_2$ initial data  is globally solved with
 \be\label{E:inverse-u-t} 
\begin{split}
|u(\cdot,t)|_{ L^1_2}=&|\mathcal S^{-1}( \{R_+(k,t),i\beta_{ 1},\cdots,i\beta_{n},\gamma_{+, 1}(t),\cdots,\gamma_{+,n}(t) \})|_{ L^1_2}\\
=&|\mp\partial B_\pm(\cdot,0^\pm,t)/\partial x|_{ L^1_2}\\
\le &Cc_{\pm}(t). 
\end{split}\ee

  But, generically  
             \be\label{E:cohn} 
             \textit{$\sup_{t\gtrless 0}c_{\pm}(t)=\infty $ }
             \ee{\cite[Lemma 1.4, 1.8,  Theorem (Agranovich \& Marchenko), Theorem (Faddeev) ]{C79}, \cite{C82}}. We have no global $L^1_2$-estimates for the Cauchy problem of the KdV equation via the IST approach. 
%%%%%%%%%%%%%%%%%%%%%%%%%%%%%%%%%%%%%%%%%%%%%%%%%%%%%%%%%%%%%%%%%%%%%%%%
\section{The deformed Gelfand-Levitan-Marchenko equation}\label{S:perturbation}
%%%%%%%%%%%%%%%%%%%%%%%%%%%%%%%%%%%%%%%%%%%%%%%%%%%%%%%%%%%%%%%%%%%%%%%%

To  prove a stability of   an $n$-soliton $u_0(x,t)$, defined by \eqref{E:n-soliton-tau}, let
\be\label{E:solitary-perturbation}
\begin{gathered}
u_v(x,0 )=u_0(x,0 )+v(x)    
\end{gathered} \ee be the initial data of a KdV solution $u_v(x,t)$, it amounts to justify continuity %of the forward scattering transform $\mathcal S(u_v(x,0))$ and continuity 
of the inverse scattering transform $\mathcal S^{-1}( \{R_{v,+}(k,t),i\beta_{v, 1},\cdots$, $i\beta_{v,\mu(v)},\gamma_{v,+, 1}(t),\cdots,\gamma_{v,+,\mu(v)}(t) \})$ in terms of the perturbation $v$. Here 
\be\label{E:sd-t}
\mathcal S(u_v(x,0))= \{R_{v,+}(k),i\beta_{v, 1},\cdots,i\beta_{v,j},\gamma_{v,+, 1},\cdots,\gamma_{v,+,j} \}.\ee For simplicity, we will abbreviate $R_{v,+}(k)$, $\gamma_{v,+,\mu(v)}$, $B_{v,+}(x,y)$, $\cdots$ by $R_{v }(k)$, $\gamma_{v, \mu(v)}$, $B_{v }(x,y),\cdots$, and $\mu(v)$ by $\mu$ if no confusion is caused from now on. As is shown in \eqref{E:inverse-u-t} and \eqref{E:cohn}, additional conditions are demanded to prove a continuity in the IST approach.

\begin{lemma}\label{L:eigenfunction} Let $u_v(x,t)$ be a KdV solution with initial data   \eqref{E:solitary-perturbation} and $ \sum_{h=0}^1|v^{(h)} |_{ L^1_2 \cap L^\infty } \ll 1$.  
  Then eigenfunctions $
       - f'' _{v,\pm}    +u_v(x,0) f_{v,\pm}   =    k^2f_{v,\pm}  $, $f_{v,\pm} = m_{v,\pm} e^{\pm ikx}\to   e^{\pm ikx}$ as $ x\to \pm\infty$, satisfy 
\be\label{E:eigen-function-continuity}
\begin{gathered}
    |  m_{ v,\pm}   -m_{0,\pm} | 
\le      \frac{C(1+\max(\mp x,0))}{1+|k| }| v|_{L^1_2},   \\
 | m'_{v,\pm}-m'_{0,\pm} |\le   C   |  v |_{   L^1_2\cap L^\infty} ,\ \
 | m''_{v,\pm}-m''_{0,\pm} |\le C \sum_{h=0}^1|v^{(h)} |_{ L^1_2 \cap L^\infty }
 \end{gathered}  
  \ee  for $x\in\RR,\,k\in\mathbf C$. 

\end{lemma}

\begin{proof} From
\begin{align}
   & m_{ v,+} -m_{ 0,+} =  {-} \int_x^{ \infty}\frac{e^{2ik  (s-x) }-1}{2ik}v(s)m_{ v,+}(s,k)ds \label{E:m-diff}     \\
&\hskip1in {-}  \int_x^{ \infty}\frac{e^{2ik (s-x)}-1}{2ik}u_0(s)(m_{v,+}(s,k)-m_{ 0,+}(s,k))ds,  \nonumber
\end{align}one obtain
\begin{align}
 &| m_{ v,+} -m_{ 0,+}| \le  \int_x^{ \infty}(s-x)|v(s)m_{ v,+}(s,k)|ds     \label{E:M-est-R}\\
&\hskip1in {+}  \int_0^{ \infty}s|u_0(s)(m_{v,+}(s,k)-m_{ 0,+}(s,k))|ds \nonumber\\
&\hskip1in {+}  \int_x^{ \infty}\max(-x,0)|u_0(s)(m_{v,+}(s,k)-m_{ 0,+}(s,k))|ds. \nonumber
\end{align}

 Applying \eqref{E:m-infty}, \eqref{E:m-diff}, and the Volterra integral argument in \cite[Lemma 2.1]{DT79}, we derive
 \be 
| m_{ v,+} -m_{ 0,+} |\le  C\frac{e^{\rho_0(x)} }{1+|k|}\int_x^\infty (s-x)(1+|s|)|v(s)|ds, \label{E:m-diff-pre}
\ee with $\frac{e^{\rho_0(x)} }{1+|k|}\int_x^\infty (s-x)(1+|s|)|v(s)|ds$ a non increasing function of $x$.

Hence \eqref{E:M-est-R} turns into
\begin{align}
&| m_{ v,+} -m_{ 0,+}| \le \frac{C(1+\max(-x,0))}{1+|k|}|v|_{L^1_2}
\label{E:M-est-R-new}\\
&\hskip1in {+}  \int_x^{ \infty}\max(-x,0)|u_0(s)(m_{v,+}(s,k)-m_{ 0,+}(s,k))|ds\nonumber
\end{align}and one derives  the estimate for $| m_{ v,+} -m_{ 0,+}|$ of \eqref{E:eigen-function-continuity}. 

Using 
\[
m'_{v,+} (x,k)=-\int_x^\infty e^{2ik(s-x)}u_v(s)m_{v,+}(s,k)ds
\] and similar argument as above, one can derive remaining estimates of \eqref{E:eigen-function-continuity}.

\end{proof}

\begin{lemma} \label{L:dsd-est}
Let $u_v(x,t)$ be a KdV solution with initial data   \eqref{E:solitary-perturbation}, there exists $C$ such that if $  |v   |_{ L^1_2 }\le\frac{\epsilon^\sigma}C$, $\sigma>2$, $ \color{black} \epsilon< \min_{1\le j\ne j'\le n}\{  \frac {\beta_{0,1}}2,\frac {|\beta_{0,j}-\beta_{0,j'}|}2 \}  $,  and $\mathfrak d $ being the greatest integer less than $1+\int_{-\infty}^\infty|s u_0(s,0)|ds$, then
\begin{align}
 n\le \mu(v)\le  &\mathfrak d,  \label{E:discrete-1}  \\
  |\beta_{v,\mu-n+j} -\beta_{0,j}|< &  \epsilon, \qquad\ \ \, 1\le j\le n, \label{E:discrete-2} \\
  |\gamma_{v,\mu-n+j} -\gamma_{0,j}|< &  C\epsilon, \qquad 1\le j\le n, \label{E:discrete-3}\\
  |\beta_{v,j'} |< &  \epsilon ,  \qquad \ \ \,1\le j'\le \mu-n . \label{E:discrete-4}
\end{align}

\end{lemma}
  
\begin{proof} Estimate \eqref{E:discrete-1} follows from 
\[
\textit{number of poles of $T_v$}\le 1 +\int_{-\infty}^\infty|s u_v|ds
\]\cite[pp.\,149]{DT79} and $| v|_{L^1_2}\le    \epsilon^\sigma/C$.

Since $i\beta_{v,1},\cdots,i\beta_{v,\mu  }$ are  simple zeros of $a_v(k)$   
\be\label{E:a-v-new}
\begin{split}
 &a_v(k) \equiv \frac 1{T_v(k)}=-\frac{1}{2ik}\textrm{Wr}\,(f_{v,+}(x, k),f_{v,-}(x, k)) \\
 = &- \frac{1}{2ik}\left(
\begin{array}{cc}
 m_{v,+}(x ,k) &   m_{v,-}(x ,k)\\
ik m_{v,+}(x ,k)+ \partial_x m_{v,+}(x ,k)&-ik m_{v,-}(x ,k)+ \partial_x m_{v,-}(x ,k)
\end{array}
\right) ,
\end{split}\ee \cite{DT79}
 and $m(x,\cdot)-1\in H^{2+}$. Hence 
 \be\label{E:a-0}
\begin{split}
 &a_0(k) =- \frac{1}{2ik}\left(
\begin{array}{cc}
 m_{0,+}(x_0,k) &   m_{0,-}(x_0,k)\\
ik m_{0,+}(x_0,k)+ \partial_x m_{0,+}(x_0,k)&-ik m_{0,-}(x_0,k)+ \partial_x m_{0,-}(x_0,k)
\end{array}
\right) ,
\end{split}\ee and $a_0\to 1$ as $|k|\to\infty$. Here we have used that $a_v(k)$ is independent of $x$. 
For convenience of estimating \eqref{E:beta-mu} later, additionally, $x_0$ is chosen to satisfy
\be\label{E:non-degenerate}
\sum_{j=1}^n \frac 1{|m_{0,+}(x_0,i\beta_{0,j})|} <C. 
\ee 

Therefore, for $ \color{black} \epsilon< \min_{j\ne j'}\{  \frac {\beta_{0,1}}2,\frac {|\beta_{0,j}-\beta_{0,j'}|}2 \}  $, {\color{black} using the simple zero property  of $\beta_{0,j}$},  there exists $\color{black}C_1$,
\be\label{E:epsilon-delta}  
\begin{gathered}
\{k:|k-\beta_{0,j}|\le \epsilon \}\cap\{k:|k-\beta_{0,j'}|\le \epsilon \}=\phi,\quad 1\le j<j'\le n,\\
|a_0(k)|>\frac \epsilon {\color{black}C_1} ,\quad\textit{if } k\in \cap _{j=1}^n \{k:|k-\beta_{0,j}|\ge \epsilon \}.%,\\\beta_{0,j}\gg  \epsilon,\quad 1\le j \le n. 
\end{gathered}
\ee 

From Lemma \ref{L:eigenfunction}, there exists $C\gg 1$, such that if $|v   |_{ L^1_2 }\le\frac{\epsilon^\sigma}C$ then
\be \label{E:zero}
 |  m_{ v,\pm}(x_0,k)  -m_{ 0,\pm} (x_0,k) |+|  \partial_x m_{ v,\pm}(x_0,k)  - \partial_x m_{ 0,\pm} (x_0,k) |  \le   {\epsilon^\sigma }   .
\ee  Along with \eqref{E:a-v-new} and enlarge $C$ if necessary, yields
\be\label{E:a-zero}
\begin{gathered}
 |  a_{ v}(k)  -a_{ 0} ( k) | 
\le    { \epsilon^{\sigma -1}}  ,\quad\textit{for }|k|\ge\epsilon. 
\end{gathered}
\ee Combining \eqref{E:epsilon-delta} and \eqref{E:a-zero}, and   using  the simple zero property of $a_0$,  we obtain, 
\begin{align}
& \textit{
if $\ \ k\in \cap _{j=1}^n \{k:|k-\beta_{0,j}|\ge \epsilon \}\cap \{|k|\ge \epsilon\}\quad $   }\label{E:a-der}\\
&\ \textit{then}\quad|a_0(k)|> \frac \epsilon {\color{black}C_1}, \ \  |  a_{ v}(k)  -  a_{ 0} ( k) |\le    \epsilon^{\sigma -1}   \textit{ and}
 \nonumber \\
  & \textit{
if $\ \ k\in \cup _{j=1}^n \{k:|k-\beta_{0,j}|\le \epsilon \}  $   }\label{E:a-der-1}\\
&\ \textit{then}\quad|  a_0(k)|\le   2 \epsilon  \max_{1\le j\le n}\sup_{|k-\beta_{0,j}|\le\epsilon} |\partial_ka_0(k)|,\nonumber\\
&\ |\partial_k a_0(k)|\le |\partial_k a_0(k)|+2\max_{1\le j\le n}\sup_{|k-\beta_{0,j}|\le\epsilon}|\frac{\partial^2}{\partial k^2}a_0(k)|, \nonumber\\
&\ |   a_{ v}(k)  -  a_{ 0} ( k) | 
\le      \frac{\epsilon^{\color{black}\sigma-1  } }{C_2}  ,\ | \partial_k a_{ v}(k)  -\partial_k a_{ 0} ( k) | 
\le      \frac{\epsilon^{\color{black}\sigma-2  } }{C_2} . \nonumber
\end{align}  

From \eqref{E:a-der}  and $\sigma>2$, we conclude
\be\label{E:zero-1} 
\textit{if $ \beta_{v,j}\ge\epsilon$ then $\beta_{v,j}\in \cup_{j'=1}^n\{k:|k-\beta_{0,j'}|\le  \epsilon \}$.}
\ee %This proves \eqref{E:discrete-4}. 

Moreover, notice that 
\eqref{E:epsilon-delta} implies that $ i\beta_{0,j} $ are the only zero (simple) of $a_0(k)$ in $\{k:|k-\beta_{0,j}|\le  \epsilon \}$. Together with   $a_{v}-1\in H^{2+}$, and a  local mapping theorem  of holomorphic functions \cite[Chapter 4, Theorem 10]{Ah53},
\be\label{E:local-map-0}
\begin{gathered}
\frac 1{2\pi i}\int_{|k-\beta_{0,j}|= \epsilon }\frac{\partial_k a_0(k)}{a_0(k)}dk=1 .
\end{gathered}
\ee
Applying \eqref{E:a-v-new},  \eqref{E:zero},  \eqref{E:local-map-0}, and the local mapping theorem,
\be \label{E:local-map-v}
  \frac 1{2\pi i}\int_{|k-\beta_{0,j}|= \epsilon}\frac{\partial_k a_v(k)}{a_v(k)}dk=\frac 1{2\pi i}\int_{|k-\beta_{0,j}|= \epsilon}\frac{\partial_k a_0(k)}{a_0(k)}dk+\mathcal O( \epsilon^{\sigma -2}  )  ,
\ee which is continuous and an integer-valued function.  Enlarge $C$, one can have $C_2$  large enough and then $|\mathcal O( \epsilon^{\sigma -2}  )|\le \epsilon^{\sigma -2}<1$. Therefore,   $\frac 1{2\pi i}\int_{|k-\beta_{0,j}|= \epsilon}\frac{\partial_k a_v(k)}{a_v(k)}dk=1$ and there is one and only one (simple) zero of $a_v $ in $\{k:|k-\beta_{0,j}|\le \epsilon\}$. Along with \eqref{E:zero-1}, we prove \eqref{E:discrete-2} and \eqref{E:discrete-4}.

Finally,  using \cite[Chapter Two, (2.2.35)]{AC91}, we obtain
\be\label{E:beta-mu}
\gamma_{v,j}=\frac {-i}{\frac d{dk}a_v(i\beta_{v,j})}\frac{m_{v,-}(x_0,i\beta_{v,j})}{m_{v,+}(x_0,i\beta_{v,j})}e^{2\beta_{v,j}x_0}.
\ee Here we have used $\gamma_{v,j}$ are constant (independent of $x$ particularly). 
Hence, for $\mu-n+1\le j\le\mu$, combining with \eqref{E:discrete-2},   \eqref{E:non-degenerate},  \eqref{E:zero}, and \eqref{E:a-der-1},  one derives    \eqref{E:discrete-3}.

\end{proof}  

The forward scattering transform $\mathcal S$ is not continuous since the number of the bound states could increase under arbitrarily small perturbation by \eqref{E:discrete-1} (see examples in \cite[Chapter 4, Section 11]{D83}).  In this case, there is no uniform positive lower bound for the bound states $\beta_{v,j}$ which cause difficulty in deriving estimates. For instance, 
        a standard Riemann-Hilbert approach involves with replacing the mixed data problem by a purely continuous problem.  Small bound states will cause the replacement failing to be a good approximation \cite{BC84}. %Besides, small bound states also form an obstruction in estimating the GLM equation and soliton dynamics (see proofs of  Lemma \ref{L:S-v} and \ref{L:soliton-phase} if).  

\begin{lemma} \label{L:csd-est}
Let $u_v(x,t)$ be a KdV solution with initial data   \eqref{E:solitary-perturbation}, and $  
 \sum_{h=0}^2|v^{(h)} e^{\mathfrak a|x|}|_{ L^\infty }$ $\le  \frac{\epsilon^\sigma}C $, $\mathfrak a>0,\,\sigma>2$, $\color{black}  \epsilon<\max\{\mathfrak a, \beta_{0,1} \}$. %, and $C$ is defined by Lemma \ref{L:dsd-est}.  
 Then $R_v$ can be   meromorphically  extended on the strip $   0\le 
k_I\le \epsilon$, and  
   \begin{eqnarray}
    R_0(k_R+ik_I)\equiv 0,\hskip1in && k_R\in\mathbf R, \quad 0\le k_I\le \epsilon,\label{E:continuous-1}\\
 | (k_R+i\epsilon)^{N+1}R_{v }(k_R+i\epsilon) |< C \sum_{h=0}^{ N}| v^{(h)}|_{L^1_2}, && k_R\in\mathbf R,\quad N=0,1,2. \label{E:continuous-2}    
\end{eqnarray}   

\end{lemma}
\begin{proof} From the assumption $ 
 \sum_{h=0}^2|v^{(h)} e^{\mathfrak a|x|}|_{ L^\infty }\le \frac{\epsilon^\sigma}C $, $\color{black}\epsilon<\max\{\mathfrak a, \beta_{0,1} \} $, \eqref{E:m-infty}, \eqref{E:alg-con}, Lemma \ref{L:eigenfunction}, and the representation formula
\be\label{E:reflection}
\begin{split}
            \frac{R_v(k)}{T_v(k)}=&   \frac{1}{2ik}\int e^{-2iks}u_v(s){m_{v,-}}(s,k)ds, \\
    \frac 1{T_v(k)}=&1 -\frac{1}{2ik}\int u_v(s){m_{v,+}}(s,k)ds,\ (T_v-1)^{-1}\in H^{2+},
        \end{split}
        \ee \cite{DT79},  one can meromorphically extend the   reflection coefficient $T_v(k)$ and 
    $ T_0(k)\equiv 1$, $ R_0(k) \equiv 0 $, on the strip $  0\le 
k_I\le \epsilon $.   %and derive the decaying estimates on $k_I=\epsilon$. 

From \eqref{E:continuous-1}, \eqref{E:reflection}, decompose
\begin{align}
  & \frac 1{2i(k_R+i\epsilon)}\int e^{-2i(k_R+i\epsilon)s}u_v(s)m_{ v,-}(s, k_R+i\epsilon)ds \label{E:R-estimste}\\
= &\frac{1}{2i(k_R+i\epsilon)}\int e^{-2i(k_R+i\epsilon)s} v(s)m_{ v,-}(s, k_R+i\epsilon)ds \nonumber\\
+& \frac{1}{2i(k_R+i\epsilon)}\int e^{-2i(k_R+i\epsilon)s} u_0(s)(m_{ v,-}(s,  k_R+i\epsilon )-m_{ 0,-}(s, k_R+i\epsilon))ds.  \nonumber
\end{align} 

For each term on the right hand side of \eqref{E:R-estimste}, one can derive estimates \eqref{E:continuous-2} by   the assumption $ 
 \sum_{h=0}^2|v^{(h)} e^{\mathfrak a|x|}|_{ L^\infty }\le \frac{\epsilon^\sigma}C $,  $\color{black}\epsilon<\max\{\mathfrak a, \beta_{0,1} \} $,  Lemma \ref{L:eigenfunction}, integration by parts,   and induction.
 
%For the first term, using the assumption    Lemma \ref{L:eigenfunction}, integration by parts,   and induction , one can derive estimates \eqref{E:continuous-2}.

\end{proof}

%Note that taking $N=0$ in the decay estimate \eqref{E:continuous-2},    we have \be\label{E:reflection-0} | R_v(k ) | \le      \min\,\{1,\,\frac C{|k|  }|v|_{L^1_2}\} \ee which is not small near $k=0$ no matter how small the perturbation $v$ is. %It causes obstruction to derive smallness estimates for a Riemann-Hilbert problem approach (see \cite[Theorem 8.9]{BC84} for instance).

From Section \ref{S:KdV-IST}, one define 
 \begin{equation} \label{E:GLM-kernel-t-new}
 \Omega_v(x,t)=  \frac 1\pi\int_{-\infty}^\infty R_v(k) e^{  2ik(x+4k^2 t)}dk+2\sum_{j=1}^\mu\gamma_{v,j}e^{- 2\beta_{v,j}(x{ - 4\beta_{v,j}^2) t}}\equiv  \Omega_{v,c}(x,t)+\Omega_{v,d}(x,t)  , 
\end{equation}  solve the GLM equation, and derive the KdV solution
    
\begin{gather}
\Omega_v(x+y,t)+B_v(x,y,t)+\int_0^{ \infty} \Omega_v (x+y+u,t)B_v(x,u,t)du=0,\ \  y\in \mathbf R^\pm ,\label{E:GLM-t-new}\\
  u_v(x,t)=-\frac{\partial}{\partial x}B_v(x,0+,t).\label{E:GLM-t-new-q}
 \end{gather}

We prove $B_v(x,y,t)$ satisfies a deformed GLM equation:  
\begin{theorem}\label{T:deforemed-GLM}  
   Let $u_v(x,t)$ be a KdV solution with initial data   \eqref{E:solitary-perturbation},   $  
 \sum_{h=0}^2|v^{(h)} e^{\mathfrak a|x|}|_{ L^\infty }$ $\le  \frac{\epsilon^\sigma}C $, $\mathfrak a>0,\,\sigma>2$, $\color{black}\epsilon < \min_{1\le j\ne j'\le n}\{\mathfrak a, \frac {\beta_{0,1}}2,\frac {|\beta_{0,j}-\beta_{0,j'}|}2 \} $, and $C$ is defined by Lemma \ref{L:dsd-est}. Hence  %the solution $B_v(x,y,t)$ of  
  \be \label{E:deform-GLM}
{ \Delta_{v  }(x+y,t)+B_{v }(x,y,t)+\int_0^\infty  \Delta_{v } (x+y+u,t)B_{v }(x,u,t)du=0,\  y>0,}
\end{equation}  with $\Delta_{v\  }(x;t)=\Delta_{v,c  }(x;t)+\Delta_{v,d  }(x;t)$  and  
 \begin{align}
 { \Delta_{v,c  }(x;t) =}& \frac { e^{-2\epsilon(x-4\epsilon^2 t)}}\pi \int_{-\infty}^\infty R_v(k+i\epsilon)e^{2ik [x+4 (k^2 -3 \epsilon^2)t] { -24 \epsilon k^2    t} }dk, \label{E:deform-GLM-c}\\
 { \Delta_{v,d  }(x;t) =}& 2\sum_ {j={ \mu-n+1}}^\mu \  \gamma_{v,j} e^{-2\beta_{v,j}(x-4\beta_{v,j}^2t)  },  \label{E:deform-GLM-d}
  \end{align}    satisfying
    \begin{align}
          \sup_{ x\ge 4\epsilon^2t{\color{black}\ge 0 }  ,y,u>0} |\Delta_{v,c}(x+y+u,t)  |_{ L^\infty  }\le &  C \epsilon^{ \color{black}\sigma -1}e^{-2\epsilon u}, \label{E:deform-GLM-c-est}\\
           %\sup_{ x\ge 4\epsilon^2t\ge 0  ,y>0} |\Delta_{v,c}(x+y+u,t)  |_{L^1(\mathbf R^+,du)  }\le &   C \epsilon^{\frac\sigma 2-1},\label{E:deform-GLM-c-est-2}\\
           \sup_{ x\ge 4\epsilon^2t{\color{black}\ge 0 }  ,y,u>0}  |\frac{\partial\Delta_{v,c}}{\partial x}(x+y+u,t)   |_{  L^\infty  }\le &  C \epsilon^{\sigma }e^{-2\epsilon u}.\label{E:deform-GLM-c-est-1}%\\
    %\sup_{ x\ge 4\epsilon^2t  ,y>0}  |\frac{\partial\Delta_{v,c}}{\partial x}(x+y+u,t)   |_{L^1(\mathbf R^+,du)  }\le &  C \epsilon^{\sigma-1}e^{-2\epsilon u}.\label{E:deform-GLM-c-est-3}
    \end{align}

\end{theorem}

\begin{proof} Assume \eqref{E:deform-GLM}-\eqref{E:deform-GLM-d} are true. Applying \corr{\eqref{E:continuous-2} to \eqref{E:deform-GLM-c}}, we derive, for $x\ge 4\epsilon^2t{\color{black}\ge 0 }$, $y,u>0$, 
\begin{align}
&|\Delta_{v,c}(x+y+u,t)|_{   L^\infty} \label{E:delta-est-1}\\
\le & C|R_v\corr{(k+i\epsilon)}|_{L^1(\RR,dk) } e^{-2\epsilon(y+u)}\nonumber\\
=&C\left[\corr{|R_v(k+i\epsilon)|_{L^1(|k|\le 1)}  }+|R_v\corr{(k+i\epsilon)}|_{L^1(|k|\ge 1)}\right] e^{-2\epsilon(y+u)}\nonumber\\
\le & C(\corr{\epsilon^{-1}}|v| _{L^1_2}+|v'| _{L^1_2}) e^{-2\epsilon(y+u)} \le C \epsilon^{ \color{black}\sigma -1 } e^{-2\epsilon(y+u)} ,\nonumber\\ 
 & |\frac{\partial}{\partial x}\Delta _{v,c}(x+y+u,t)|_{   L^\infty }\label{E:delta-est-2}\\
 \le & C|kR_v|_{L^1(\RR,dk) } e^{-2\epsilon(y+u)}+C \epsilon^{\frac\sigma 2+1 } e^{-2\epsilon(y+u)} \nonumber\\
 \le & C\sum_{h=0}^2|v^{(h)}|_{L^1_2} e^{-2\epsilon(y+u)}+C \epsilon^{\frac\sigma 2+1 } e^{-2\epsilon(y+u)} \nonumber\\
 \le & C \epsilon^{\sigma } e^{-2\epsilon(y+u)} . \nonumber 
\end{align}
Therefore follow \eqref{E:deform-GLM-c-est} and \eqref{E:deform-GLM-c-est-1}. It then reduces to proving \eqref{E:deform-GLM}-\eqref{E:deform-GLM-d}. 
To this aim, we introduce the notation  and properties
\be\label{E:fourier}
\begin{split}
\hat f(y)={f\,}^{\widehat{•}}=\frac 1\pi\int_{-\infty}^\infty f(k)e^{2iky}dk,\qquad& 
\check g(k)= \int_{-\infty}^\infty g(y)e^{-2iky}dy,\\
\check{\hat f}(k)=f(k),\qquad&  \hat{\check g}(y)= g(y),\\
\widehat{f_1 f_2} = \hat f_1\ast\hat f_2 ,\qquad & \widehat{f_1\ast f_2} =\frac 1\pi\hat f_1 \hat f_2 , \\
\frac 1\pi\int_{-\infty}^\infty f(k)\overline{g(k)}dk =&\int_{-\infty}^\infty \hat f(y)\overline{\hat g}(y)dy,
\end{split}\ee and 
\be\label{E:fourier-1}
\begin{split}
\gamma_{v,j}(0)=\left(\int_{-\infty}^\infty f_{v,+}^2(x,i\beta_{v,j} )dx\right)^{-1},\qquad & f_{v,-}(x,t,i\beta_{v,j})=c_{v,j}(t)f_{v,+}(x,t,i\beta_{v,j}),\\
R_v(k,t)= R_v(k)e^{8ik^3 t},\qquad &T_v(k,t)= T_v(k,0),\\
\gamma_{v,j}(t)= \gamma_{v,j}(0)e^{8\beta_{v,j}^3 t},\qquad & c_{v,j}(t)= c_{v,j}(0)e^{8\beta_{v.j }^3t} .
\end{split}\ee

Taking the Fourier transform of both  sides of 
\be\label{E:system-extend}  {T_v (k,t )}f_{v,-}(x,t,k )=  {R_{v,+} (k,t )} f_{v,+}(x,t,k )+ f_{v,+}(x,t,- k ),\ee  applying Lemma \ref{L:eigenfunction}, \ref{L:csd-est}, and using the residue theorem, we derive
\begin{align}
&\frac 1\pi\int_{-\infty}^\infty(T_v(k+i\epsilon,t)-1)f_{v,-}(x,t,k+i\epsilon)e^{2i(k+i\epsilon)y}dk\label{E:sd-KdV}\\
=&\frac 1\pi\int_{-\infty}^\infty R_v(k+i\epsilon,t)(f_{v,+}(x,t,k+i\epsilon)-e^{ i(k+i\epsilon)x})e^{2i(k+i\epsilon)y}dk\nonumber\\
+&\frac 1\pi\int_{-\infty}^\infty R_v(k+i\epsilon,t) e^{ i(k+i\epsilon)(x+2y)}dk\nonumber\\
+&\frac 1\pi\int_{-\infty}^\infty (f_{v,+}(x,t,-(k+i\epsilon))-e^{ -i(k+i\epsilon)x})e^{2i(k+i\epsilon)y}dk\nonumber\\
-&\frac 1\pi\int_{-\infty}^\infty (f_{v,-}(x,t, k+i\epsilon )-e^{ -i(k+i\epsilon)x})e^{2i(k+i\epsilon)y}dk.\nonumber
\end{align}

Therefore, using \cite[pp. 55 (24)]{T74}, \cite[pp. 146]{DT79},  Lemma \ref{L:dsd-est}, the residue theorem,   \eqref{E:fourier}, and \eqref{E:fourier-1},
\begin{align}
 \blacktriangleright\quad&\frac 1\pi\int_{-\infty}^\infty(T_v(k+i\epsilon,t)-1)f_{v,-}(x,t,k+i\epsilon)e^{2i(k+i\epsilon)y}dk\label{E:sd-KdV-LHS}\\
=&2i\sum_{j=\mu-n+1}^\mu \res_{k=i\beta_{v,j}}  (T_v(k)-1)f_{v,-}(x,t,k)e^{2iky} \nonumber\\
=&2i\sum_{j=\mu-n+1}^\mu i\left(\int_{-\infty}^\infty f_{v,+}(x,i\beta_{v,j})f_{v,-}(x,i\beta_{v,j})dx\right)^{-1}    f_{v,-}(x,t,i\beta_{v,j})e^{-2\beta_{v,j}y} \nonumber\\
=&-2\sum_{j=\mu-n+1}^\mu  \gamma_{v,j}(0)c_{v,j}^{-1}(0)     f_{v,-}(x,t,i\beta_{v,j})e^{-2\beta_{v,j}y}\nonumber \\
 =&-2\sum_{j=\mu-n+1}^\mu  \gamma_{v,j}(t)c_{v,j}^{-1}(t)     f_{v,-}(x,t,i\beta_{v,j})e^{-2\beta_{v,j}y}\nonumber \\
=&-2\sum_{j=\mu-n+1}^\mu  \gamma_{v,j}(t)      f_{v,+}(x,t,i\beta_{v,j})e^{-2\beta_{v,j}y} \nonumber\\
=&-2\sum_{j=\mu-n+1}^\mu  \gamma_{v,j}(t)\left(e^{-\beta_{v,j}x}+ e^{-\beta_{v,j}x} \int_{0}^\infty B_v(x,\tilde y,t) e^{-2\beta_{v,j}\tilde y}d\tilde y\right)e^{-2\beta_{v,j}y}\nonumber\\
=&-2\sum_{j=\mu-n+1}^\mu  \gamma_{v,j}(t)\left(e^{-2\beta_{v,j}(y+\frac x2)}+  \int_{0}^\infty B_v(x,\tilde y,t) e^{-2\beta_{v,j}(\tilde y+y+\frac{x}{2})}d\tilde y\right) ;\nonumber\\
&\nonumber\\
 \blacktriangleright\quad&\frac 1\pi\int_{-\infty}^\infty R_v(k+i\epsilon,t)(f_{v,+}(x,t,k+i\epsilon)-e^{ i(k+i\epsilon)x})e^{2i(k+i\epsilon)y}dk\label{E:sd-KdV-RHS-1} \\
=&\int_{-\infty}^\infty \left(R_v(k+i\epsilon,t)e^{2i(k+i\epsilon)y}\right)^{\widehat {}}(u)\  \left(f_{v,+}(x,t,k+i\epsilon)-e^{ i(k+i\epsilon)x}\right)^{\widehat{ }  }(-u)dk\nonumber\\
=&\int_{-\infty}^\infty \left(R_v(k+i\epsilon)e^{8i(k+i\epsilon)^3t}e^{2i(k+i\epsilon)y}\right)^{\widehat {}}(u)\  \left(f_{v,+}(x,t,k+i\epsilon)-e^{ i(k+i\epsilon)x}\right)^{\widehat{ }  }(-u)du\nonumber \nonumber\\
=& \int_{-\infty}^\infty  e^{-2\epsilon y}\left(R_v(k+i\epsilon)e^{8ik (k^2 -3 \epsilon^2)t +8 \epsilon( \epsilon^2  -3k^2 ) t }\right)^{\widehat{•}}(u+y)
\nonumber\\
\times &\left[\frac{1}{\pi}\int_{-\infty}^\infty \left(
e^{i(k+i\epsilon)x}\int_0^\infty B_v(x,\tilde y,t)e^{2i(k+i\epsilon)\tilde y}d\tilde y  \right)e^{-2iku}dk\right]du \nonumber\\
=& \int_{-\infty}^\infty   e^{-2\epsilon (u+y)} \left(R_v(k+i\epsilon)e^{8ik  (k^2 -3 \epsilon^2)t +8 \epsilon( \epsilon^2   -3k^2   )t }\right)^{\widehat{•}}(u+y)
   B_v(x,u-\frac{x}{2},t) du; \nonumber\\
   &\nonumber\\
    \blacktriangleright\quad&\frac 1\pi\int_{-\infty}^\infty R_v(k+i\epsilon,t)e^{2i(k+i\epsilon)(y+\frac{x}{2})}dk\label{E:sd-KdV-RHS-2} \\
%=&\frac{1}{\pi}\int_{-\infty}^\infty  R_v(k+i\epsilon)e^{8i(k+i\epsilon)^3t}  e^{2i(k+i\epsilon)(y+\frac{x}{2})}dk \nonumber\\
%=&\frac{1}{\pi}\int_{-\infty}^\infty  R_v(k+i\epsilon)e^{8i(k ^3t-3k\epsilon^2 t+k(y+\frac x2)}  e^{-(24 k^2\epsilon t -8\epsilon^3 t+2\epsilon (y+\frac{x}{2})}dk \nonumber  \\
=& e^{-2\epsilon (y+\frac{x}{2})}\left(R_v(k+i\epsilon)e^{8ik  (k^2 -3 \epsilon^2)t +8 \epsilon( \epsilon^2   -3k^2   )t }\right)^{\widehat{•}}(y+\frac x2);\nonumber\\&\nonumber\\
   \blacktriangleright\quad&\frac 1\pi\int_{-\infty}^\infty (f_{v,+}(x,t,-(k+i\epsilon))-e^{ -i(k+i\epsilon)x})e^{2i(k+i\epsilon)y}dk\label{E:sd-KdV-RHS-3} \\
  =& \frac 1\pi\int_{-\infty}^\infty \left(e^{ -i(k+i\epsilon)x}\int_0^\infty B_v(x,\tilde y,t) e^{ -2i (k+i\epsilon)\tilde y}d\tilde y\right)e^{2i(k+i\epsilon)y}dk\nonumber\\
  =&   B_v(x,  y-\frac x2,t)  ;\nonumber\\&\nonumber\\
  \blacktriangleright\quad&-\frac 1\pi\int_{-\infty}^\infty (f_{v,-}(x,t, k+i\epsilon )-e^{ -i(k+i\epsilon)x})e^{2i(k+i\epsilon)y}dk\label{E:sd-KdV-RHS-4}  \\
  =& -\frac 1\pi\int_{-\infty}^\infty \left(e^{ -i(k+i\epsilon)x}\int_{-\infty}^0    B_{v,-}(x,\tilde y,t) e^{ -2i (k+i\epsilon)\tilde y}d\tilde y\right)e^{2i(k+i\epsilon)y}dk\nonumber\\
  =& -   B_{v,-}(x,  y-\frac x2,t)  .\nonumber
\end{align}

Combining \eqref{E:sd-KdV}-\eqref{E:sd-KdV-RHS-4}, we obtain
\begin{align*}
&-2\sum_{j=\mu-n+1}^\mu  \gamma_{v,j}(t) e^{-2\beta_{v,j}(y+\frac x2)}
-2\sum_{j=\mu-n+1}^\mu  \gamma_{v,j}(t) \int_{0}^\infty B_v(x,\tilde y,t) e^{-2\beta_{v,j}(\tilde y+y+\frac{x}{2})}d\tilde y  \\
=&\int_{-\infty}^\infty   e^{-2\epsilon (u+y)} \left(R_v(k+i\epsilon)e^{8ik  (k^2 -3 \epsilon^2)t +8 \epsilon( \epsilon^2   -3k^2   )t }\right)^{\widehat{•}}(u+y)
   B_v(x,u-\frac{x}{2},t) du\\
   +&e^{-2\epsilon (y+\frac{x}{2})}\left(R_v(k+i\epsilon)e^{8ik  (k^2 -3 \epsilon^2)t +8 \epsilon( \epsilon^2   -3k^2   )t }\right)^{\widehat{•}}(y+\frac x2)\\
   +&  B_v(x,  y-\frac x2,t)- B_{v,-}(x,  y-\frac x2,t).
\end{align*}
Setting $
y-\frac{x}{2}=z>0$ on both sides, $
u-\frac x2=\tilde z>0$ on RHS, $\tilde y=\tilde z>0$ on LHS, the above identity turns into
\begin{align*}
&-2\sum_{j=\mu-n+1}^\mu  \gamma_{v,j}(t) e^{-2\beta_{v,j}(z+x)}
-2\sum_{j=\mu-n+1}^\mu  \gamma_{v,j}(t) \int_{0}^\infty B_v(x,\tilde z,t) e^{-2\beta_{v,j} (z+\tilde z+x)}d\tilde z  \\
=&\int_{-\infty}^\infty   e^{-2\epsilon (z+\tilde z+x)} \left(R_v(k+i\epsilon)e^{8ik  (k^2 -3 \epsilon^2)t +8 \epsilon( \epsilon^2   -3k^2   )t }\right)^{\widehat{•}}(z+\tilde z+x)
   B_v(x,\tilde z,t) d\tilde z\\
   +&e^{-2\epsilon (z+x)}\left(R_v(k+i\epsilon)e^{8ik  (k^2 -3 \epsilon^2)t +8 \epsilon( \epsilon^2   -3k^2   )t }\right)^{\widehat{•}}(z+x) 
   +   B_v(x,  z,t).  
\end{align*} Therefore we complete the proof of \eqref{E:deform-GLM}-\eqref{E:deform-GLM-d}.
\end{proof}

To summarize, advantages to replace the GLM equation \eqref{E:GLM-t-new} by the deformed GLM equation \eqref{E:deform-GLM} include: instead of \eqref{E:cohn} of $\Omega_{v,c}$ in \eqref{E:GLM-kernel-t-new}, we gain the {integrability} of $\Delta_{v,c}$ in a moving regime,  as well as, despite the forward scattering transform is discontinuous,  $\Delta_{v,d}$ is free of extra small bound states of $\Omega_{v,d}$ and   is close to the GLM kernel $\Omega_{0}$ of the unperturbed $n$-soliton $u_0(x,t)$ due to Lemma \ref{L:dsd-est}.

%%%%%%%%%%%%%%%%%%%%%%%%%%%%%%%%%%%%%%%%%%%%%%%%%%%%%%%%%%%%%%%%%%%%%%%%%
\section{Proof of Theorem \ref{T:stability}}\label{S:stability}
%%%%%%%%%%%%%%%%%%%%%%%%%%%%%%%%%%%%%%%%%%%%%%%%%%%%%%%%%%%%%%%%%%%%%%%%%

%In this section, several necessary conditions of the IST are provided for solving the inverse problem in the next section.

To derive a stability from the deformed GLM equation, we define
\be\label{E:GLM-operator}
\begin{split}
 \mathcal K_{v,c}f =&  \int_0^\infty \Delta_{v,c}(x+y+u,t) f(x,u,t) du, \\
\mathcal K_{v,d}f =&\int_0^\infty \Delta_{v,d}(x+y+u,t) f(x,u,t) du,
\end{split}
\ee and
\be\label{E:GLM-solution}
\begin{split}
    (1+\mathcal K_{v,d})B_{v,d}=&-\Delta_{v,d},\\
    B_{v,c}=&B_{v}-B_{v,d} .
\end{split}
\ee %Note {$\mp\frac{\partial }{\partial x}B_{v,d}(x,0^\pm)$ are KdV $n$-solitons with bdd states $i\beta_{v,\mu-n+1},\cdots,i\beta_{v,\mu}$ and norming constants $\gamma_{v,\mu-n+1},\cdots,\gamma_{v,\mu}$}.
Here  $\mp\frac{\partial }{\partial x}B_{v,d}(x,0^\pm)$ are KdV $n$-solitons with bound states $i\beta_{v,\mu-n+1},\cdots,i\beta_{v,\mu}$ and norming constants $\gamma_{v,\mu-n+1}$, $\cdots,\gamma_{v,\mu}$.  The deformed GLM equation  \eqref{E:deform-GLM} formally turns into
\begin{multline}     (1+ (1+\mathcal K_{v,d})^{-1}\mathcal K_{v,c})B_{v,c}\\ = -(1+ (1+\mathcal K_{v,d})^{-1}\mathcal K_{v,c}) B_{v,d} -(1+\mathcal K_{v,d})^{-1}(\Delta_{v,d} +\Delta_{v,c}) .\label{E:GLM-decom}\end{multline}
Therefore, we are led to prove invertibility of $1+\mathcal K_{v,d}$, derive estimates of $ \mathcal K_{v,c}$, $(1+\mathcal K_{v,d})^{-1}$, $B_{v,d}$, and $B_{v,c}$.

\begin{lemma}\label{L:K-c} Let $u_v(x,t)$ be a KdV solution with initial data   \eqref{E:solitary-perturbation},  $  
 \sum_{h=0}^2|v^{(h)} e^{\mathfrak a|x|}|_{ L^\infty }$ $\le  \frac{\epsilon^\sigma}C $, $\mathfrak a>0,\,\sigma>2$, $\color{black}\epsilon < \min_{1\le j\ne j'\le n}\{\mathfrak a, \frac {\beta_{0,1}}2,\frac {|\beta_{0,j}-\beta_{0,j'}|}2 \} $, and $C$ is defined by Lemma \ref{L:dsd-est}.  Then 
\[\begin{split}
| \mathcal K_{v,c}f|_{L^\infty (x\ge 4\epsilon^2t{\color{black}\ge 0 } ,y>0)} \le   C \epsilon^{\color{black} \sigma -2}|f|_{L^\infty (x\ge 4\epsilon^2t {\color{black}\ge 0 },y>0)},\\
 |\mathcal K' _{v,c}f|_{L^\infty (x\ge 4\epsilon^2t{\color{black}\ge 0 } ,y>0)}  \le C \epsilon^{\sigma-1}|f|_{L^\infty (x\ge 4\epsilon^2t {\color{black}\ge 0 },y>0)}. 
\end{split}\]
\end{lemma}
\begin{proof} From \eqref{E:delta-est-1} and \eqref{E:delta-est-2},
\begin{align*}
\blacktriangleright\quad&| \mathcal K_{v,c}f|_{L^\infty (x\ge 4\epsilon^2t{\color{black}\ge 0 } ,y>0)}\\
=&| \int_0^\infty \Delta_{v,c}(x+y+u,t) f(x,u,t) du|_{L^\infty (x\ge 4\epsilon^2t{\color{black}\ge 0 } ,y>0)}\\
\le&|f|_{L^\infty (x\ge 4\epsilon^2t {\color{black}\ge 0 },y>0)}  \int_0^\infty |\Delta_{v,c}(x+y+u,t) |_{L^\infty (x\ge 4\epsilon^2t {\color{black}\ge 0 },y,u>0)} du\\
%\le &C|f|_{L^\infty (x\ge 4\epsilon^2t ,y>0)}| \int_0^\infty\epsilon^{\frac\sigma 2 } e^{-2\epsilon(y+u)} du|_{L^\infty ( y>0)}\\
\le &C|f|_{L^\infty (x\ge 4\epsilon^2t{\color{black}\ge 0 } ,y>0)}  \int_0^\infty \epsilon^{ \color{black}\sigma -1 } e^{-2\epsilon u} du \\
\le & C \epsilon^{\color{black} \sigma -2}|f|_{L^\infty (x\ge 4\epsilon^2t{\color{black}\ge 0 } ,y>0)}, \\
\blacktriangleright\quad&| \mathcal K'_{v,c}f|_{L^\infty (x\ge 4\epsilon^2t{\color{black}\ge 0 } ,y>0)}\\
=&| \int_0^\infty \Delta'_{v,c}(x+y+u,t) f(x,u,t) du|_{L^\infty (x\ge 4\epsilon^2t{\color{black}\ge 0 } ,y>0)}\\
\le&|f|_{L^\infty (x\ge 4\epsilon^2t{\color{black}\ge 0 } ,y>0)}  \int_0^\infty |\Delta'_{v,c}(x+y+u,t) |_{L^\infty (x\ge 4\epsilon^2t {\color{black}\ge 0 },y,u>0)} du\\
%\le &C|f|_{L^\infty (x\ge 4\epsilon^2t ,y>0)}| \int_0^\infty\epsilon^{\frac\sigma 2 } e^{-2\epsilon(y+u)} du|_{L^\infty ( y>0)}\\
\le &C|f|_{L^\infty (x\ge 4\epsilon^2t{\color{black}\ge 0 } ,y>0)}  \int_0^\infty \epsilon^{ \sigma   } e^{-2\epsilon u} du \\
\le & C \epsilon^{ \sigma  -1}|f|_{L^\infty (x\ge 4\epsilon^2t {\color{black}\ge 0 },y>0)}.
\end{align*} 
\end{proof}

The following lemma   is an extension  of \cite[Lemma 6.1]{ES83} where they justified, for  $c$, $M>0$ fixed, there exists $C_{u_0,c,M}$ such that   $1+\mathcal K_{0,d}$ is invertible on $L^\infty(\mathbf R \times\mathbf R^+\times\mathbf R )$ and $ S_0=(1+\mathcal K_{0,d})^{-1}$ satisfies
  \be\label{E:es-est} 
| S_0 f|_{ L^\infty(| x- 4c^2 t|\le M,x, y,t\ge 0)},\,| S' _0f|_{ L^\infty(| x- 4c^2 t|\le M, x, y,t\ge 0)}\le C_{u_0,c,M}|f|_{ L^\infty(| x- 4c^2 t|\le M, x, y,t\ge 0)}.
\ee The estimates \eqref{E:es-est} are important to yield an asymptotic theorem by taking $t\to\infty$ along lines of fixed slopes  for any fixed KdV solution via IST approach \cite{ES83}.

The following lemma shows that   uniform boundedness of $S_v=(1+\mathcal K_{v,d})^{-1}$ can be justified in the soliton region with cones around the $x$-axis and lines determined by bound states of  $u_0(x,0)$ removed. Therefore, we can eventually conclude an orbital stability theorem with explicit phase shifts in a regime where $t$ is not necessarily large.

\begin{lemma}\label{L:S-v}

Given positive constants   $\mathfrak a$, $ \sigma>2$,   $  \color{black}\tau $, there exists a constant $C\gg 1$ such that if $u_v(x,t)$  a KdV solution with initial data   \eqref{E:solitary-perturbation}, for $\color{black}\epsilon < \min_{1\le j\ne j'\le n}\{\mathfrak a, \frac {\beta_{0,1}}2,\frac {|\beta_{0,j}-\beta_{0,j'}|}2 \} $, then    the transform $1+\mathcal K_{v,d}$ is invertible on   $\color{black}L^\infty( \frac x{4t}\in \cap_{j=1}^n \{c^2:|c-\beta_{0,j}|\ge \tau,  \frac 1\tau\ge c \},t,y\ge 0)$ and $ S_v=(1+\mathcal K_{v,d})^{-1}$ satisfying%. Writing , 
\begin{align} 
  | S_v f|_{\color{black}L^\infty( \frac x{4t}\in \cap_{j=1}^n \{c^2:|c-\beta_{0,j}|\ge \tau,  \frac 1\tau\ge c\},t,y\ge 0)} \le &C|f|_{\color{black}L^\infty( \frac x{4t}\in \cap_{j=1}^n \{c^2:|c-\beta_{0,j}|\ge \tau,  \frac 1\tau\ge c\},t,y\ge 0)},\label{E:S-v}\\
 |   S' _v f|_{\color{black}L^\infty( \frac x{4t}\in \cap_{j=1}^n \{c^2:|c-\beta_{0,j}|\ge \tau,  \frac 1\tau\ge c\},t,y\ge 0)}\le& C|f|_{\color{black}L^\infty( \frac x{4t}\in \cap_{j=1}^n \{c^2:|c-\beta_{0,j}|\ge \tau,  \frac 1\tau\ge c\},t,y\ge 0)}. \label{E:S-v-1}
\end{align} %for $ \sum_{h=0}^2|v^{(h)} e^{\mathfrak a|x|}|_{ L^\infty }   \le  \frac{\epsilon^\sigma}C $, $\mathfrak a>0,\,\sigma>2$, and  $\color{black}\epsilon < \epsilon_0 $.

\end{lemma}  
\begin{proof}Write
  \begin{align}
   \Delta_{v,d}(x+y+u,t)= &   \sum_{j=\mu-n+1}^\mu  L_{v,j}(x+u,t) e^{-2\beta_{v,j} y   }, \nonumber\\ 
 L_{v,j}(x+u,t)= & 2\gamma_{v,j} e^{-2\beta_{v,j}(x+u-4\beta_{v,j}^2t)  },\nonumber\\
   \mathfrak  L_{v,j} = &L_{v,j}(x,t) \nonumber\\
   =&\corr{2\gamma_{v,j} e^{-2\beta_{v,j}(x-4c^2t)+4(c^2-\beta_{v,j}^2)t)  }}\nonumber\\
   = &\corr{ G_j(x-4c^2 t) e^{-8\beta_{v,j}(c^2-\beta_{v,j}^2)t}},\nonumber\\
   \mathcal K_{v,d}f= &  \sum_{j=\mu-n+1}^\mu A_{v,j}(x,t)e^{-2\beta_{v,j} y   }, \label{E:k-vd}\\ 
 A_{v,j}(x,t)=&  \mathfrak L_{v,j} \int_0^\infty e^{-2\beta_{v,j} u   }f(x,u,t)du.\nonumber
\end{align}

Following argument in \cite{ES83}, one can prove $I+\mathcal K_{v,d}$ is invertible by  
  \be\label{E:s-v} 
S_vf =( I+\mathcal K_{v,d})^{-1}f =f(x,y,t)-\sum_{j= 1}^n A_{v,\mu-n+j}(x,t)e^{-2\beta_{v,\mu-n+jj} y   },
\ee and
\be\label{E:a-v} 
\begin{gathered}
   \Gamma_v  \left(
\begin{array}{c}
A_{v, \mu-n+1}\\
A_{v, \mu-n+2}\\
\vdots \\
A_{v, \mu }
\end{array}
\right)= 
\left(
\begin{array}{c}
2\int_0^\infty e^{-2\beta_{v,\mu-n+1} u   }f(x,u,t)du\\
2\int_0^\infty e^{-2\beta_{v,\mu-n+2} u   }f(x,u,t)du\\
\vdots \\
2\int_0^\infty e^{-2\beta_{v,\mu} u   }f(x,u,t)du
\end{array}
\right), 
\end{gathered}
\ee with
\be\label{E:gamma-v} 
\begin{gathered}
\Gamma_v=   \left(
\begin{array}{cccc}
\frac{2}{\mathfrak L_{v,\mu-n+1}}+\frac{1}{2\beta_{v,\mu-n+1}} &\frac{1}{\beta_{v,\mu-n+1}+\beta_{v,\mu-n+2}}&\cdots & \frac{1}{\beta_{v,\mu-n+1}+\beta_{v,\mu }}\\
\frac{1}{\beta_{v,\mu-n+1}+\beta_{v,\mu-n+2}}&\frac{2}{\mathfrak L_{v,\mu-n+2}}+\frac{1}{2\beta_{v,\mu-n+2}} &\cdots & \frac{1}{\beta_{v,\mu-n+2}+\beta_{v,\mu }}\\
\vdots&\vdots&\ddots&\vdots\\
\frac{1}{\beta_{v,\mu-n+1}+\beta_{v,\mu }} &\frac{1}{\beta_{v,\mu-n+2}+\beta_{v,\mu }}&\cdots & \frac{2}{\mathfrak L_{v,\mu }}+\frac{1}{2\beta_{v,\mu }}
\end{array}
\right)  
\end{gathered} \ee an invertible matrix.

Applying Cramer's rule, 
\be\label{E:cramer} 
|A_{v,j}|\le \frac {|f|_{L^\infty(\mathbf R^+)}}{\beta_{v,\mu-n+1}}\sum_{i=1}^n\left| \frac{\widetilde \Gamma_{v,ij}}{\mbox{det}\,\Gamma_v} \right|.
\ee Here $\widetilde\Gamma_{v,ij}$ denotes the cofactor of $(i,j)$ entry in $\Gamma_v$ which is a symmetric matrix.

{\color{black}Given $v(x)$, $M$, $c$, %applying Lemma \ref{L:dsd-est} and 
following argument in \cite{ES83}, one can prove that 
\be\label{E:l-v-j}
\lim_{t\to\infty}\mathfrak L_{v,\mu-n+j} 
=\left\{{\begin{array}{ll}
\infty &\textit{if $c<\beta_{v,\mu-n+j}$}\\
G_{v,\mu-n+j} &\textit{if $c=\beta_{v,\mu-n+j}$}\\
0&\textit{if $c>\beta_{v,\mu-n+j}$}
\end{array}}
\right.
\ee for $|x-4c^2t|\le M, t>0$. By writing
\be\label{E:cramer-new-1-new}
\dfrac{\widetilde\Gamma _{v,ij}}{\mbox{det}\,\Gamma^l_v}=
\left\{ 
{\begin{array}{ll}
 \dfrac{\widetilde\Gamma _{v,ij}}{\mbox{det}\,\Gamma _v}  ,&\,\textit{ if $0<c\le\beta_{v,\mu-n+1}$},\\
 %&\\
  \dfrac{\mathfrak L _{v,\mu-n+1}\cdots\mathfrak  L _{v,\mu-n+m}\widetilde\Gamma_{v,ij} }{\mathfrak L _{v,\mu-n+1}\cdots\mathfrak  L _{v,\mu-n+m}\mbox{det}\,\Gamma_v } ,&{\begin{array}{l}\textit{if $\beta_{v,\mu-n+m}<c\le \beta_{v,\mu-n+m+1}$,  } \\
    m=1,\cdots,n-1,\end{array}}\\
 \dfrac{\mathfrak L _{v,\mu-n+1}\cdots \mathfrak L _{v,\mu}\widetilde\Gamma _{v,ij}}{\mathfrak L _{v,\mu-n+1}\cdots \mathfrak L _{v,\mu}\mbox{det}\,\Gamma _v} ,&\, \textit{ if $\beta_{v,\mu}<c$,}
\end{array}}
\right.
\ee   then   
\begin{itemize}
\item [$\blacktriangleright$] for $v$, $M$, $c$ fixed,  $\dfrac{\widetilde\Gamma _{v,ij}}{\mbox{det}\,\Gamma _v}$ converge uniformly as $t\to\infty$ on $|x-4c^2t|\le M$;
\item  [$\blacktriangleright$]  for different $c$,    the above limits (explicit formula are given in \cite{ES83} and  convergent speeds are different.
\end{itemize}
  
  Furthermore, applying Lemma \ref{L:dsd-est}, 
 \begin{itemize}
\item [$\blacktriangleright$]   $\dfrac{\widetilde\Gamma _{v,ij}}{\mbox{det}\,\Gamma _v}$ converge uniformly as $t\to\infty$  on $\frac x{4t}\in \cap_{j=1}^n \{c^2:|c-\beta_{0,j}|\ge \tau \}$  
\end{itemize}for $\forall v$ satisfying $\sum_{h=0}^2|v^{(h)} e^{\mathfrak a|x|}|_{ L^\infty }\le     \frac{\epsilon^\sigma}C $, $ \epsilon < \min_{  j\ne j' }\{\mathfrak a, \frac {\beta_{0,1}}2,\frac {|\beta_{0,j}-\beta_{0,j'}|}2 \}$. Namely, there exists constant $A$, such that 
\be\label{E:unif-new}
\begin{split}
 &\textit{$|\dfrac{\widetilde\Gamma _{v,ij}}{\mbox{det}\,\Gamma _v}|<C$, for $t>A$, $\frac x{4t}\in \cap_{j=1}^n \{c^2:|c-\beta_{0,j}|\ge \tau \}$ } 
\end{split}
\ee for $\forall v$ satisfying $\sum_{h=0}^2|v^{(h)} e^{\mathfrak a|x|}|_{ L^\infty }\le     \frac{\epsilon^\sigma}C $, $ \epsilon < \min_{  j\ne j' }\{\mathfrak a, \frac {\beta_{0,1}}2,\frac {|\beta_{0,j}-\beta_{0,j'}|}2 \}$.

Since   $\frac 1\tau\ge \sqrt{\frac x{4t}}$ and $0\le t\le A$ yield $x\le \frac{4A}{\tau^2}$.
Together with the invertibility of $\Gamma_0$, Lemma \ref{L:dsd-est}, and enlarging $C$ if necessary, gives rise to,  
\be\label{E:unif-new-1}
\begin{split}
&\textit{$|\dfrac{\widetilde\Gamma  _{v,ij}}{\mbox{det}\,\Gamma  _v}|<C$, for      $0\le t \le A,\,  x\le \frac{4A}{\tau^2} $  } 
\end{split}
\ee for $\forall v$ satisfying $\sum_{h=0}^2|v^{(h)} e^{\mathfrak a|x|}|_{ L^\infty }\le     \frac{\epsilon^\sigma}C $, $ \epsilon < \min_{  j\ne j' }\{\mathfrak a, \frac {\beta_{0,1}}2,\frac {|\beta_{0,j}-\beta_{0,j'}|}2 \}$.

Combining \eqref{E:unif-new} and \eqref{E:unif-new-1}, for $\forall v$ satisfying $\sum_{h=0}^2|v^{(h)} e^{\mathfrak a|x|}|_{ L^\infty }\le     \frac{\epsilon^\sigma}C ,\, \epsilon < \min_{  j\ne j' }\{\mathfrak a$, $ \frac {\beta_{0,1}}2$, $\frac {|\beta_{0,j}-\beta_{0,j'}|}2 \}$, 
\be\label{E:s-v-00-new} 
\begin{gathered}
|\dfrac{\widetilde \Gamma _{v,ij}}{\mbox{det}\,\Gamma _v}|_{L^\infty(\frac x{4t}\in \cap_{j=1}^n \{c^2:|c-\beta_{0,j}|\ge  \tau ,\frac 1\tau\ge c\},t,y\ge 0  )} <C,
\end{gathered}\ee
 and   
\be\label{E:s-v-0}
\begin{gathered}
|A_{v,j}|_{L^\infty(\frac x{4t}\in \cap_{j=1}^n \{c^2:|c-\beta_{0,j}|\ge  \tau ,\frac 1\tau\ge c\},t,y\ge 0  )} <C  |f|_{L^\infty(\frac x{4t}\in \cap_{j=1}^n \{c^2:|c-\beta_{0,j}|\ge  \tau ,\frac 1\tau\ge c\},t,y\ge 0  )},\\
| S_vf|_{L^\infty(\frac x{4t}\in \cap_{j=1}^n \{c^2:|c-\beta_{0,j}|\ge  \tau ,\frac 1\tau\ge c\},t,y\ge 0  )}\le C  |f|_{L^\infty(\frac x{4t}\in \cap_{j=1}^n \{c^2:|c-\beta_{0,j}|\ge  \tau ,\frac 1\tau\ge c\},t,y\ge 0  )}.
\end{gathered}
\ee}

The proof of \eqref{E:S-v-1} follows from the identities \cite[(6.21),  (6.22)]{ES83}
  \[
  \begin{split}
   S_v'f  =&-\sum_{j=1}^n A'_{v,\mu-n+j}(x,t) e^{-2\beta_{v,\nu-n+j}y},\\ 
    A'_{v,\mu-n+j}   =& \sum_{i=1}^n   \frac{-4\beta_{v,\mu-n+i} A_{v,\mu-n+i} }{\mathfrak L_{v,\mu-n+i}}\dfrac{\widetilde \Gamma _{v,ij}}{\mbox{det}\,\Gamma _v},
  \end{split}
  \]and similar argument as above.   
\end{proof}

\begin{lemma} \label{L:estimate-b} Given positive constants   $\mathfrak a$, $ \sigma>2$,   $  \color{black}\tau $, there exists $C\gg 1$ such that if $u_v(x,t)$ is a KdV solution satisfying \eqref{E:solitary-perturbation},   $  
 \sum_{h=0}^2|v^{(h)} e^{\mathfrak a|x|}|_{ L^\infty } \le  \frac{\epsilon^\sigma}C $,   $ \epsilon < \min_{  j\ne j' }\{\mathfrak a, \frac {\beta_{0,1}}2,\frac {|\beta_{0,j}-\beta_{0,j'}|}2 \}$, then
{ \begin{align}
  |B_{v,d}|_{\color{black}L^\infty(\frac x{4t}\in \cap_{j=1}^n \{c^2:|c-\beta_{0,j}|\ge  \tau ,\frac 1\tau\ge c\},t,y\ge 0  )}\le & C,\label{E:b-d-0}\\
     | B' _{v,d}|_{\color{black}L^\infty(\frac x{4t}\in \cap_{j=1}^n \{c^2:|c-\beta_{0,j}|\ge  \tau ,\frac 1\tau\ge c\},t,y\ge 0  )}\le &  C,\label{E:b-d}\\
 |B_{v,c}|_{\color{black}L^\infty(\frac x{4t}\in \cap_{j=1}^n \{c^2:|c-\beta_{0,j}|\ge  \tau ,\frac 1\tau\ge c\ge\epsilon\},t,y\ge 0  )}\le& C\epsilon^{\color{black} \sigma -2 }\label{E:b-c}\\
   | B '_{v,c}|_{\color{black}L^\infty(\frac x{4t}\in \cap_{j=1}^n \{c^2:|c-\beta_{0,j}|\ge  \tau ,\frac 1\tau\ge c\ge\epsilon\},t,y\ge 0  )}\le  &  C\epsilon^{\color{black} \sigma -2}. \label{E:b-c-'}
\end{align}}
As a result,
\begin{equation}\label{E:dev} 
\begin{gathered}
%\sup_{ x\ge 4\epsilon^2 t }|B_v-B_{v,d}|_{L^\infty(\mathbf R^+,dy)} \le  C\epsilon ^{\sigma-1},\\
 |u_v-B'_{v,d}(x,0^+,t) |_{\color{black}L^\infty(\frac x{4t}\in \cap_{j=1}^n \{c^2:|c-\beta_{0,j}|\ge  \tau ,\frac 1\tau\ge c\ge\epsilon\},t > 0  )}\le  C\epsilon ^{\color{black} \sigma -2}.
\end{gathered}
\end{equation}
\end{lemma}
 
\begin{proof} Estimates for $\sup_{  (x,t)\in(\mathbf R^+,\mathbf R^+)} |B_{v,d}|_{L^\infty(\mathbf R^+,dy)}  $, $ \sup_{  (x,t)\in(\mathbf R^+,\mathbf R^+)} | B'_{v,d}|_{L^\infty(\mathbf R^+,dy)} $ follow from Lemma \ref{L:dsd-est}, the identities \cite[(6.26), (6.28), (6.31)]{ES83}
  \[
  \begin{gathered}
    B_{v,d}(x,y,t)= \sum_{j=1}^n (A_{v,\nu-n+j}(x,t)-L_{v,\nu-n+j}(x,t))e^{-2\beta_{v,\nu-n+j}y},\\ 
    A_{v,\nu-n+j}(x,t)-L_{v,\nu-n+j}(x,t) = -2\sum_{i=1}^n \frac{\Gamma_{v,ij}}{\det \Gamma_{v}} ;\end{gathered}
  \]
    and
    \[
  \begin{gathered}
     B'_{v,d}(x,y,t)= \sum_{j=1}^n ( A'_{v,\nu-n+j}(x,t)- L'_{v,\nu-n+j}(x,t))e^{-2\beta_{v,\nu-n+j}y},\\ 
 A'_{v,\nu-n+j}(x,t)- L'_{v,\nu-n+j}(x,t)=  \sum_{j=1}^n \frac{-4\beta_{v,\mu-n+i}(A_{v,\nu-n+i}-L_{v,\nu-n+i})}{L_{v,\mu-n+i}}\frac{\Gamma_{v,ij}}{\det \Gamma_{v}},
 \end{gathered}
  \]and similar argument as above.

  Estimates for $ B_{v,c} $, $  B '_{v,c}$ follow from  applying Lemma \ref{L:K-c}, \ref{L:S-v}, \eqref{E:GLM-solution}, and \eqref{E:b-d}   to  \eqref{E:GLM-decom} and its $x$-derivative.

   Finally, \eqref{E:dev} follows from \eqref{E:GLM-t-new-q}, \eqref{E:GLM-solution},  and \eqref{E:b-c-'}.

\end{proof}

The following technical lemma describes dynamics of KdV $n$-solitons under perturbation. We prove it via an application of the asymptotic theory of KdV $n$-solitons \cite{K04, K18}. 
\begin{lemma}\label{L:soliton-phase}
Given positive constants   $\mathfrak a$, $ \sigma>2$,  there exists $C\gg 1$ such that if $u_v $ is a KdV solution satisfying \eqref{E:solitary-perturbation}, %{\color{black} there exists $\epsilon_2  \ll \min\{\mathfrak a, \beta_{0,1} \}$} such that 
  $  
 \sum_{h=0}^2|v^{(h)} e^{\mathfrak a|x|}|_{ L^\infty } \le  \frac{\epsilon^\sigma}C $ for  
 $ \color{black}\epsilon < \min_{  j\ne j' }\{\mathfrak a, \frac {\beta_{0,1}}2,\frac {|\beta_{0,j}-\beta_{0,j'}|}2 \} $, then  
\begin{multline}\label{E:soliton-est}
 \quad|2\frac{d^2}{dx^2}\log \textrm{Wr}\left(a_{  1}(x-x_{ 1}(t),t),\cdots,a_{  n}(x-x_{ n}(t),t)\right)\\
  -2\frac{d^2}{dx^2}\log \textrm{Wr}\left(a_{v,  1}(x ,t),\cdots,a_{v,  n}(x ,t)\right)| \le C\epsilon,\quad
 \end{multline}
where
\be\label{E:phase-final}
 x_{ j}(t)=4 (\beta^2_{v,\mu-n+j} -\beta^2_{0,j})t,\quad j=1,\cdots,n,
 \ee and 
\begin{gather} 
 a_{v,j}(x,t)=(-1)^{j+1}e^{-\beta_{v,\mu-n+j} (x-4\beta_{v,\mu-n+j}^2 t)}+\alpha_{v,\mu-n+j}e^{+\beta_{v,\mu-n+j} (x-4\beta_{v,\mu-n+j}^2 t)}, \label{E:Eckhaus-n-soliton-a_j-new} \\
  \alpha_{v,j}= \frac{2(-1)^{j-1}\beta_{v,j}}{\gamma_{v,j}}\left(\Pi_{l\ne j}\frac{\beta_{v,l}+\beta_{v,j}}{\beta_{v,l}-\beta_{v,j}}\right),\quad 1\le j\le \mu(v), 
 \label{E:Eckhaus-n-soliton-alpha-new}\\
 a_j\equiv a_{0,j},\quad \alpha_{ j}\equiv\alpha_{0,j}, \quad 1\le j\le \mu(0)=n,\label{E:new-v-0}\\
 \textit{$\beta_{v,j}$, $\gamma_{v,j}$    are bound states, norming constants  of $u_v(x,0)$}.\nonumber
\end{gather}

\end{lemma}
\begin{proof} Abbreviating  $\textrm{Wr}(f_1,\cdots,f_n)=\left|\begin{array}{ccc} f_1&\cdots&f_n \\\vdots&\ddots&\vdots\\
\frac{d^{(n-1)}}{dx^{n-1}}f_1&\cdots &\frac{d^{(n-1)}}{dx^{n-1}}f_n\end{array}\right| $ by $\left|\begin{array}{c} f_j\\\vdots\\\frac{d^{(n-1)}}{dx^{n-1}}f_j\end{array}\right|_{j=1,\cdots,n}$,  let's decompose 
\begin{multline}\label{E:soliton-est-1}
\qquad 2\frac{d^2}{dx^2}\log \textrm{Wr}\left(a_{  1}(x-x_{ 1}(t),t),\cdots,a_{  n}(x-x_{ n}(t),t)\right)\\
  -2\frac{d^2}{dx^2}\log \textrm{Wr}\left(a_{v,  1}(x ,t),\cdots,a_{v,  n}(x ,t)\right)=I_{v,1}+I_{v,2}, \qquad
 \end{multline} with 
 \begin{align}
 &I_{v,1}\label{E:I-1}\\
 = & \textit{\footnotesize  $2\frac{d^2}{dx^2}\log \left|
\begin{array}{r}
(-1)^{j+1}e^{-\beta_j \{x-4\beta_j^2 t+4(\beta_j^2-\beta_{v,\mu-n+j}^2)t\}}+\alpha_je^{+\beta_j \{x-4\beta_j^2 t+4(\beta_j^2-\beta_{v,\mu-n+j}^2)t\}}\\
%\frac{d}{dx}[(-1)^{j+1}e^{-\beta_j \{x-4\beta_j^2 t+4(\beta_j^2-\beta_{v,\mu-n+j}^2)t\}}+\alpha_je^{+\beta_j \{x-4\beta_j^2 t+4(\beta_j^2-\beta_{v,\mu-n+j}^2)t\}}]\\
\vdots \hskip2.05 in \\
\frac{d^{(n-1)}}{dx^{n-1}}[(-1)^{j+1}e^{-\beta_j \{x-4\beta_j^2 t+4(\beta_j^2-\beta_{v,\mu-n+j}^2)t\}}+\alpha_je^{+\beta_j \{x-4\beta_j^2 t+4(\beta_j^2-\beta_{v,\mu-n+j}^2)t\}}]
\end{array}
 \right|_{j=1,\cdots,n}$} \nonumber\\
 &\textit{\footnotesize $-2\frac{d^2}{dx^2}\log \left|
\begin{array}{r}
(-1)^{j+1}e^{-\beta_{v,\mu-n+j }(x-4\beta_{v,\mu-n+j }^2 t)}+\alpha_{j }e^{+\beta_{v,\mu-n+j } (x-4\beta_{v,\mu-n+j }^2 t)}\\
%\frac{d}{dx}[(-1)^{j+1}e^{-\beta_{v,\mu-n+j } (x-4\beta_{v,\mu-n+j }^2 t)}+\alpha_{j }e^{+\beta_{v,\mu-n+j } (x-4\beta_{v,\mu-n+j }^2 t)}]\\
\vdots \hskip2.16in \\
\frac{d^{(n-1)}}{dx^{n-1}}[(-1)^{j+1}e^{-\beta_{v,\mu-n+j } (x-4\beta_{v,\mu-n+j }^2 t)}+\alpha_{j }e^{+\beta_{v,\mu-n+j } (x-4\beta_{v,\mu-n+j }^2 t)}]
\end{array}
 \right|_{j=1,\cdots,n},$} \nonumber 
 \end{align}  and
 \begin{align}
&I_{v,2}\label{E:I-2}\\
=&\textit{\footnotesize $ 2\frac{d^2}{dx^2}\log \left|
\begin{array}{r}
(-1)^{j+1}e^{-\beta_{v,\mu-n+j }(x-4\beta_{v,\mu-n+j }^2 t)}+\alpha_{j }e^{+\beta_{v,\mu-n+j } (x-4\beta_{v,\mu-n+j }^2 t)}\\
%\frac{d}{dx}(-1)^{j+1}e^{-\beta_{v,\mu-n+j } (x-4\beta_{v,\mu-n+j }^2 t)}+\alpha_{j }e^{+\beta_{v,\mu-n+j } (x-4\beta_{v,\mu-n+j }^2 t)}\\
\vdots \hskip2.16in \\
\frac{d^{(n-1)}} {dx^{n-1}}[(-1)^{j+1}e^{-\beta_{v,\mu-n+j } (x-4\beta_{v,\mu-n+j }^2 t)}+\alpha_{j }e^{+\beta_{v,\mu-n+j } (x-4\beta_{v,\mu-n+j }^2 t)}]
\end{array}
 \right|_{j=1,\cdots,n}$}\nonumber\\
 &\textit{\footnotesize $-2\frac{d^2}{dx^2}\log \left|
\begin{array}{r}
(-1)^{j+1}e^{-\beta_{v,\mu-n+j }(x-4\beta_{v,\mu-n+j }^2 t)}+\alpha_{v,\mu-n+j }e^{+\beta_{v,\mu-n+j } (x-4\beta_{v,\mu-n+j }^2 t)}\\
%\frac{d}{dx}(-1)^{j+1}e^{-\beta_{v,\mu-n+j } (x-4\beta_{v,\mu-n+j }^2 t)}+\alpha_{v,\mu-n+j }e^{+\beta_{v,\mu-n+j } (x-4\beta_{v,\mu-n+j }^2 t)}\\
\vdots \hskip2.16in \\
\frac{d^{(n-1)}}{dx^{n-1}}[(-1)^{j+1}e^{-\beta_{v,\mu-n+j } (x-4\beta_{v,\mu-n+j }^2 t)}+\alpha_{v,\mu-n+j }e^{+\beta_{v,\mu-n+j } (x-4\beta_{v,\mu-n+j }^2 t)}]
\end{array} 
 \right|_{j=1,\cdots,n}.$}\nonumber 
\end{align}  

For $n=1$, one has
\begin{align}
I_{v,1} =&\frac{\beta_1^2}2\textrm{sech}^2\frac 12(\beta_1[x-4 \beta_{v,\mu}^2t]+\ln\alpha_1)-\frac{\beta_{\mu,1}^2}2\textrm{sech}^2\frac 12(\beta_{\mu,1}[x-4\beta_{v,\mu}^2t]+\ln\alpha_1),\label{E:I_1-n-1}\\
I_{v,2} =&\frac{\beta_{\mu,1}^2}2\textrm{sech}^2\frac 12(\beta_{\mu,1}[x-4\beta_{v,\mu}^2t]+\ln\alpha_1)-\frac{\beta_{\mu,1}^2}2\textrm{sech}^2\frac 12(\beta_{\mu,1}[x-4\beta_{v,\mu}^2t]+\ln\alpha_{v,\mu}).\label{E:I_1-n-2}
\end{align}

If $
 \sum_{h=0}^2|v^{(h)} e^{\mathfrak a|x|}|_{ L^\infty } \le  \frac{\epsilon^\sigma}C $, $\mathfrak a>0,\,\sigma>2$, and 
 $ \color{black}\epsilon < \min_{1\le j\ne j'\le n}\{\mathfrak a, \frac {\beta_{0,1}}2,\frac {|\beta_{0,j}-\beta_{0,j'}|}2 \} $, $C\gg 1$, applying Lemma \ref{L:dsd-est}, there exists $N\gg 1$ (independent of $v$), such that   
\be\label{E:I_1-n-1-G}
|I_{v,1}|,\,|I_{v,2}|\le \epsilon, \ \textit{ for $ |x-4 \beta_{v,\mu}^2t|\ge N$. }
\ee 

From Lemma \ref{L:dsd-est} and compactness,  there exists $C=C_N$, such that
\be\label{E:I_1-n-1-l}
|I_{v,1}|,\,|I_{v,2}|\le C \epsilon  , \ \textit{ for $ |x-4 \beta_{v,\mu}^2t|\le N$.}
\ee 
Therefore, the lemma is proved if $n=1$.

Cases for $n>1$ can be justified in the same spirit once we apply and adapt  the asymptotic theory for KP multi line solitons \cite[Theorem 6.1, Property 7.1]{K18}. We sketch the argument and leave details to \cite{K18}.%Since proofs for them are identical, we only give the proof for $I_2$ for simplicity.

Firstly we identify KdV $n$-solitons \eqref{E:n-soliton-tau} as the KP $n$-solitons of P-type with $y\equiv 0$. Precisely, we identify the Wronskians \eqref{E:n-soliton} and \eqref{E:n-soliton-a_j}  as    
\be\label{E:tau-asy}
\begin{split}
\tau(x,t)=&\left(
\begin{array}{cccccccc}
1 &0 & \cdots & 0&0  & \cdots &0 & b_1\\
0 &1 & \cdots & 0&0 & \cdots & b_2 &0\\
\vdots & \ddots &\ddots &\vdots&\vdots & \ddots &\ddots &\vdots\\
0 &\cdots & 0 & 1&b_n & 0& \cdots  &0
\end{array}
\right)
\left(
\begin{array}{ccc}
E_1 & \cdots & \kappa_1^{n-1}E_1\\
E_2 &  \cdots  &\kappa_2^{n-1}E_2\\
\vdots &  \ddots &\vdots\\
E_{2n} &\cdots &  \kappa_{2n}^{n-1}E_{2n}
\end{array}
\right),\\
=&\sum_{1\le j_1< \cdots< j_n\le 2n}\Delta_{j_1,\cdots,j_n}(A)E_{j_1,\cdots,j_n}(x,t),\nonumber
\end{split}
\ee where
\be\label{E:line-grassmannian}
 \begin{gathered} 
 A=\left(
\begin{array}{cccccccc}
1 &0 & \cdots & 0&0  & \cdots &0 & b_1\\
0 &1 & \cdots & 0&0 & \cdots & b_2 &0\\
\vdots & \ddots &\ddots &\vdots&\vdots & \ddots &\ddots &\vdots\\
0 &\cdots & 0 & 1&b_n & 0& \cdots  &0
\end{array}
\right), \quad   b_j =(-1)^{n-j}\alpha_j,\\
E_j (x,t) =\exp\theta_{j}(x,t) ,\quad \theta_j(x,t)=\kappa_jx-4\kappa_j^3t ,\\ 
  \kappa_1=-\beta_n,\ \kappa_2=-\beta_{n-1},\ \cdots,\ \kappa_n=-\beta_1,\\  \kappa_{2n}= \beta_n,\ \kappa_{2n-1}= \beta_{n-1},\ \cdots,\ \kappa_{n+1}= \beta_1,
 \end{gathered} 
 \ee $\Delta_{j_1,\cdots,j_n}(A)$ is the $n\times n$ minor of the matrix $A$ whose columns are labelled by the   index set $J=\{j_1< \cdots< j_n\}\subset\{1,\cdots,2n\}$, and $
 E_{J}(x,t)=E_{j_1,\cdots,j_n}(x,t)=\Pi_{l<m}(\kappa_{j_m}- \kappa_{j_l})\exp ( \sum_{h=1}^n\theta_{j_h}(x,t)  )$. 
 
Therefore, replacing the $y$-asymptotic theory for KP multi line solitons \cite[Theorem 6.1, Property 7.1]{K18} by the $t$-asymptotic theory \cite{K04}, due to the tau function around the soliton is given by the dominant exponential terms
\be\label{E:tau-dominant}
\tau(x,t)\sim |\Delta_I|E_I+|\Delta_J|E_J \quad\textit{as $t\to\pm\infty$},
\ee with $ I\backslash \{j\} = J\backslash \{2n-j+1\} $, 
 the soliton has wave crest (peak) along the line solitons
\be\label{E:wave-crest}
\begin{gathered}
2\frac{d^2}{dx^2} \log \tau(x,t)\sim 2\beta_j^2\textrm{sech}^2\frac 12(\beta_j[x-4 \beta_{j}^2t]+\rho^\pm_j),\quad\textit{as $t\to\pm\infty$},\\
\rho^\pm=\log\frac{|\Pi_{l<m}(\kappa_{i_m}- \kappa_{i_l})||\Delta_I|}{|\Pi_{l<m}(\kappa_{j_m}- \kappa_{j_l})||\Delta_J|} .
\end{gathered}
\ee

 Consequently, if $
 \sum_{h=0}^2|v^{(h)} e^{\mathfrak a|x|}|_{ L^\infty } \le  \frac{\epsilon^\sigma}C $, $\mathfrak a>0,\,\sigma>2$, $C\gg 1$, and 
 $  \epsilon < \min_{1\le j\ne j'\le n}\{\mathfrak a$, $\frac {\beta_{0,1}}2,\frac {|\beta_{0,j}-\beta_{0,j'}|}2 \} $, applying \eqref{E:wave-crest} and Lemma \ref{L:dsd-est}, there exists $A\gg 1,\,N\gg 1$ (independent of $v$), such that  if   $ |x|+|t|\ge A$, 
\be\label{E:I_1-n-1-G-new}
|I_{v,1}|,\,|I_{v,2}|\le \epsilon, \ \textit{ for $ \cap_{1\le j\le n}|x-4\beta_{v,\mu-n+j}^2t|\ge N$; }
\ee 
and thanks to compactness,  there exists $C=C_N$, such that
\be\label{E:I_1-n-1-l-new-1}
|I_{v,1}|,\,|I_{v,2}|\le C \epsilon  , \ \textit{ for $\cup_{1\le j\le n}|x-4\beta_{v,\mu-n+j}^2t|\le N$.}
\ee 

Moreover, from Lemma \ref{L:dsd-est}  and compactness,  there exists $C=C_A$, such that
\be\label{E:I_1-n-1-l-new}
|I_{v,1}|,\,|I_{v,2}|\le C \epsilon  , \ \textit{ for $ |x|+|t|\le A$.}
\ee 

Therefore, the lemma is proved.
\end{proof} 

\noindent $\underline{\textbf{Proof of Theorem \ref{T:stability}}}:$
 
From \eqref{E:GLM-solution} and \cite[Theorem 3.6, Remark 3.3]{DT79}, 
\be\label{E:Eckhaus-n-soliton-a_j-new-0} 
 -\frac{\partial}{\partial x}B_{v,d}(x,0+;t)
= -2\frac{d^2}{dx^2}\log \textrm{Wr}(a_{v,1}(x,t),\cdots,a_{v,n}(x,t)).\ee

%Applying \eqref{E:discrete-2}, \eqref{E:discrete-3}, \eqref{E:Eckhaus-n-soliton-a_j-new-0}-\eqref{E:Eckhaus-n-soliton-alpha-new},   \eqref{E:dev}, and denoting where by applying Lemma \ref{L:dsd-est}  and looking at dominant $e^{\Theta_J}$, $\Theta_J=\pm\theta_{v,\mu-n+1}\cdots \pm\theta_{v,\mu}$  on the $x,t$-plane where $\theta_{v,\mu-n+j}=\beta_{v,\mu-n+j}(x-4\beta_{v,\mu-n+j}^2)t$ \cite{K18}.  

 \corr{Therefore,    on $\frac x{4t}\in \cap_{j=1}^n \{c^2:|c-\beta_{0,j}|\ge  \tau ,\frac 1\tau\ge c\ge\epsilon\} $, $t>0$,   $  \epsilon < \min_{1\le j\ne j'\le n}\{\mathfrak a,\frac {\beta_{0,1}}2$, $\frac {|\beta_{0,j}-\beta_{0,j'}|}2 \} $, applying \eqref{E:dev}  and Lemma \ref{L:soliton-phase},}
\begin{align*}
 & |u_v(x,t)+2\frac{d^2}{dx^2}\log \textrm{Wr}\left(a_{  1}(x-x_{ 1}(t),t),\cdots,a_{  n}(x-x_{ n}(t),t)\right)| \\
 \le & |u_v(x,t)+2\frac{d^2}{dx^2}\log \textrm{Wr}\left(a_{v,  1}(x ,t),\cdots,a_{v,  n}(x ,t)\right)| \\
 +&|2\frac{d^2}{dx^2}\log \textrm{Wr}\left(a_{  1}(x-x_{ 1}(t),t),\cdots,a_{  n}(x-x_{ n}(t),t)\right)\\
 &-2\frac{d^2}{dx^2}\log \textrm{Wr}\left(a_{v,  1}(x ,t),\cdots,a_{v,  n}(x ,t)\right)| \\
 \le &C\epsilon^{\color{black}  \sigma- 2  }+C\epsilon \\
 \le &C\epsilon^{\color{black}\min\{1, \sigma- 2 \}}.  
\end{align*}

\end{document}